\documentclass[12pt]{amsart}
\usepackage{amsfonts,amssymb,amsmath,amsthm}
\usepackage{cite}
\usepackage{dptsymb}
\usepackage[dvips, final]{graphics}
\usepackage{epsfig}
\usepackage{isolatin1}
\usepackage{subfigure}
\def\printname#1{
	\if\draft y
		\smash{\makebox[0pt]{\hspace{-0.5in}
			\raisebox{8pt}{\tt\tiny #1}}}
	\fi
}

\newcommand{\eepic}[2]{}
\newcommand{\silenteepic}[2]{}

\newread\testin
\def\maybeinput#1{
\openin\testin=#1
\ifeof\testin\typeout{Warning: input #1 not found}\else\input#1\fi
\closein\testin
}

\newwrite\transout
\openout\transout=figs.list
\newcommand{\figscaled}[2]{%
\write\transout{-m #2 draws/#1.fig}%
\maybeinput{draws/#1.tex}%
}

\newcommand{\figcent}[2]{\mathcenter{\figscaled{#1}{#2}}}
\newcommand{\kont}[1]{\figcent{#1}{0.4}}

\usepackage{calc}
\def\mathcenter#1{
	\raisebox{.5ex+\depth-.5\totalheight}{\hbox{#1}}
}

\newtheorem{Definition}{Definition}
\newtheorem{Lemma}{Lemma}

\newtheorem{Example}{Example}
\newtheorem{Theorem}{Theorem}
\newtheorem{Proposition}{Proposition}
\newtheorem{Corollary}{Corollary}

\newcommand{\rtree}[1]{\figcent{#1}{0.4}}
\newcommand{\sub}{\operatorname{sub}}
\newcommand{\Star}{\operatorname{star}}
\newcommand{\graph}{\operatorname{graph}}
\newtheorem{Th}{Theorem}

\theoremstyle{remark}
\newtheorem*{Ack}{Acknowledgment}

\title{\bf Formal symplectic groupoid}

\author[A.~S.~Cattaneo]{Alberto~S.~Cattaneo}
\address{Institut f\"ur Mathematik, Universit\"at Z\"urich--Irchel,
Winterthurerstrasse 190, CH-8057 Z\"urich, Switzerland}
\email{asc@math.unizh.ch}

\author[B.~Dherin]{Benoit Dherin}
\address{D-MATH, ETH-Zentrum, CH-8092 Z\"urich, Switzerland}
\email{dherin@math.ethz.ch}

\author[G.~Felder]{Giovanni Felder}
\address{D-MATH, ETH-Zentrum, CH-8092 Z\"urich, Switzerland}
\email{felder@math.ethz.ch}

\thanks{A.~S.~C. acknowledges partial support of SNF Grant No.~20-100029/1}
\thanks{B. D. and G. F. acknowledge partial support of SNF Grant No.~21-65213.01}

\begin{document}

   \maketitle

\begin{abstract}
The multiplicative structure of the trivial symplectic groupoid over $\mathbb R^d$ associated
to the zero Poisson structure can be expressed in terms of a generating function. We address the 
problem of deforming such a generating function in the direction of a non-trivial
Poisson structure so that the multiplication remains associative. We prove that
such a deformation is unique under some reasonable conditions and we give
the explicit formula for it. This formula turns out to be the semi-classical
approximation of Kontsevich's deformation formula. For the case of a
linear Poisson structure, the deformed generating function reduces exactly
to the CBH formula of the associated Lie algebra. The methods used to prove
 existence are interesting in their own right as they come from an at first
sight unrelated domain of mathematics: the Runge--Kutta theory of the numeric integration
of ODE's.
\end{abstract}


\section{Introduction}

In this paper we give a formal version of the integration of Poisson manifolds by symplectic groupoids.
The solution of this formal integration problem relies on the existence of a generating function for 
which we give here the explicit formula.
This generating function turns out to be a universal Campbell--Baker--Hausdorff(CBH) formula for the non linear case. 
It reduces to the usual CBH formula when the Poisson structure comes from a Lie algebra.
This generating function can be interpreted as the semi-classical part of
Kontsevich deformation quantization formula.
This fact reminds the origin of symplectic groupoids which were first
introduced by Weinstein in \cite{CDW1987}, Karasev in \cite{karasev1987},
and Zakrwewski in \cite{zakrzewski1990} as a tool to quantize the algebra of functions on a Poisson manifold. 
This section is devoted to recall some basic features of the program of quantization by
symplectic groupoid, to formulate
the formal integration problem for Poisson manifolds and to state the main theorem of this article which gives
a positive answer to the formal integration problem.

\subsection{Quantization by symplectic groupoid}

The program of quantization by symplectic groupoid is  an attempt
to quantize the algebra of functions on a Poisson manifolds by geometric means.

It is based mainly on the belief or hope, coming from geometric quantization,
that there should exist a kind of correspondence or dictionary between the world of symplectic
manifold (classical level) and the world of linear spaces(quantum level). 
This correspondence, as explained in \cite{weinstein1997}, is summarized in the following table:

\vspace{1cm}

\begin{tabular}{|c|c|}
\hline
Symplectic world & Linear world\\
\hline
$M$              &    $Q(M)$  \\
$L\subset M $    &    $Q(L) \in Q(M)$\\
$\overline{M}$        &    $Q(\overline{M})= Q(M)^*$\\
$Q(M\times N)$   &    $Q(M)\otimes Q(N)$\\
\hline
\end{tabular}

\vspace{1cm}

Here $M$ is a symplectic manifold,
$\overline{M}$ the same manifold with opposite symplectic structure,
$L$ a Lagrangian submanifold,
and $Q(M)$ a complex vector space. $Q$ stands for the ``Quantization functor".
In particular, canonical relations, i.e., Lagrangian
submanifolds of $\overline M\times N$ are sent by $Q$ to linear maps from
$Q(M)$ to $Q(N)$. The main ingredient is the assumption that quantization is functorial
, i.e.,   the composition of canonical relations should be sent to the composition of linear maps
(see \cite{weinstein1991}).
If such a quantization functor existed, we could ask the following question: 
 \par
\emph{To what kind of symplectic
manifold should we associate an algebra( i.e., a vector space with an associative product)?}  
\par
Answering this question leads directly to the notion of symplectic groupoid, see \cite{weinstein1996}.

\begin{Definition}
A symplectic groupoid  is a Lie groupoid $G$ (see \cite{weinstein1997} for a precise
definition of a Lie groupoid) with a symplectic form $\omega$ for which
the multiplication space $G^{(m)} =\{(x,y,x\bullet y)/x,y\in G \textrm{ are composable elements}\}$
is a Lagrangian submanifold of $\overline G\times\overline G\times G$ ($\overline G$ being the symplectic 
manifold with symplectic form $-\omega$).
It can be shown (see \cite{weinstein1987}) that, given a symplectic groupoid $G$, there is an induced Poisson structure on
the base space $G^{(0)}$.
Conversely, given a Poisson manifold $P$ we call symplectic groupoid over $P$ any symplectic groupoid
$G$ such that the base space $G^{(0)}$ is diffeomorphic as Poisson manifold to $P$. In this
case we say that $G$ integrates $P$ and we call integrable Poisson manifolds the
Poisson manifolds for which we can find such a $G$.
\end{Definition}

Applying the ``Quantization functor" $Q$ to the symplectic groupoid $G$, we should then get
a vector space $Q(G)$ and an associative product $Q(G^{(m)})$ on it. The associativity of this
product being guaranteed by the associativity of the groupoid multiplication and the functoriality
of $Q$.

These  facts suggest the following procedure to quantize  Poisson manifolds $P$:

\begin{enumerate}
\item[\textbf{Step 1}] Find a symplectic groupoid $G$ such that the base $G^{(0)}$ is diffeomorphic to the Poisson
manifold $P$ .
\item[\textbf{Step 2}] Quantize (geometric quantization,...) $G$ and $G^{(m)}$ to get the quantum algebra.
\end{enumerate}

This is the idea of quantization by symplectic groupoid. \textbf{Step 1} is known as the \textbf{integrability problem} 
and was recently completely settled. Coste, Dazord and Weinstein in \cite{CDW1987} and independently Karasev in 
\cite{karasev1987} showed 
the existence of a local symplectic groupoid over any Poisson manifold, 
``local'' meaning that the multiplication is defined only on a neighborhood
of the unit space. Cattaneo and Felder in \cite{feldcat2000} gave an explicit construction
of a topological groupoid canonically associated to any Poisson manifold, which is a global
symplectic groupoid whenever the Poisson structure is integrable.
Crainic and Fernandes 
in \cite{CF2001}  derived an if and only if criterium which tells one when the previous
construction yields a manifold. 
\textbf{Step 2} however was only partially achieved (see \cite{WX1991}).

\par
If we compare this program with deformation quantization (see \cite{FFLS1997} and \cite{kontsevich1997}), we see that starting with an integrable Poisson manifold
$P$ whose symplectic groupoid is $G$ we should have the following relation between objects involved in these programs:

\vspace{1cm}
\begin{tabular}{|c|c|c|}
                     \hline
                     & Deformation quantization                  & Quantization by symplectic groupoids\\
\hline
Semi-classical level & ?                                     & $(G,G^{(m)})$\\
Quantum level       &$(C^{\infty}(P)[[\epsilon]],*_\epsilon)$  & $(Q(G),Q(G^{(m)}))$\\
\hline
\end{tabular}
\vspace{1cm}

We can regard  the symplectic groupoid over a Poisson manifold as a (semi-)classical version of the quantum algebra.
In this picture $G^{(m)}$ should then correspond to a semi-classical version of the Kontsevich star-product
formula. 
This is in some sense the case. Namely we can restate the integrability problem into a formal integration
problem. The solution of this problem is called the \textbf{formal symplectic groupoid} over a 
Poisson manifold which is a formal version of the ``true symplectic groupoid" that exists however
even for non-integrable Poisson structures.
This is exactly what the question marks stand for in the above table. Let us be more precise.

\subsection{Formal integration problem for Poisson manifolds}

In the sequel we will only consider Poisson structures $\alpha$ over $M=\mathbb R^d$. 
Suppose that $(M,\alpha)$ is  integrable and that its symplectic groupoid $G$ satisfies the
two following properties (which are always satisfied in a neighborhood of $M$):

\begin{enumerate}
\item $ G\subset T^*M\simeq \mathbb R^{*d} \times \mathbb R^d$
\item $G^{(m)}\subset \overline{T^*M}\times \overline{T^*M}\times  T^*M$ is an exact Lagrangian manifold, i.e., there
exists a generating function $S:\mathbb R^{*d}\times \mathbb R^{*d}\times \mathbb R^d\rightarrow \mathbb R$
such that $G^{(m)}= \graph(dS)$.
\end{enumerate}

We would like to see what sort of constraints the associativity of the groupoid product imposes on $S$.
First of all we may remark that under the previous assumptions the product space $G^{(m)}$ can be described
as follows
$$
G^{(m)} = \Big\{\Big(
\big(  p_1,\nabla_{p_1}S\big),
\big(  p_2,\nabla_{p_2}S\big),
\big(  \nabla_x S,x     \big)
\Big):(p_1,p_2,x)\in B_2\Big\}
$$
where the partial derivatives are evaluated at $(p_1,p_2,x)\in B_2:= (\mathbb R^{*d})^2\times \mathbb R^d$.
\par
The groupoid product associativity could be expressed by saying that, whenever the composition is allowed, we have
$g = \bar g\bullet g_3$ and $g=g_1\bullet \tilde g$ 
where $\bar g =  g_1\bullet g_2 $ and $\tilde g = g_2\bullet g_3$.

 Denoting $g=(p,x)$, $\bar g = (\bar p,\bar x)$ and $\tilde g =
(\tilde p,\tilde x)$ implies that  $(g_1,g_2,\bar g)\in G^{(m)}$, $(g_2,g_3,\tilde g)\in G^{(m)}$,
$(\bar g,g_3, g)\in G^{(m)}$ and $(g_1,\tilde g, g)\in G^{(m)}$. Now expressing $g_1,g_2,g_3,g,\bar g$ and $\tilde g$ 
each time in terms of the generating function $S$ and equating the different expressions found for the same element we get a system of six equations which can be summarized into the following more compact equation.

\vspace{1cm}

\textbf{Symplectic Groupoid Associativity equation (SGA equation):}
\par

$$S(p_1,p_2,\bar x) + S(\bar p,p_3,x)-\bar x\bar p = S(p_2,p_3,\tilde x)+S(p_1,\tilde p,x)-\tilde x\tilde p, $$

where 
$$\bar x =\nabla_{p_1}S(\bar p ,p_2,x) ,\quad \bar p =\nabla_{x}S(p_1 ,p_2,\bar x), $$
$$\tilde x =\nabla_{p_2}S(p_1,\tilde p,x) ,\quad \tilde p =\nabla_{x}S(p_2 ,p_3,\tilde x) .$$

\vspace{1cm}

This equation encodes the associativity of the groupoid product into the generating function. 
It can also be seen  from two other different points of view. First it is easy to check that one gets the SGA equation
by requiring that the saddle point evaluation  as $h$ goes to $0$ of the two integrals
$$\int e^{\frac{i}{h}[S(p_1,p_2,x)+S(p,p_3,x)-px]}\frac{d^dpd^dx}{(2\pi h)^{d/2}} \quad\textrm{ and }\quad
\int e^{\frac{i}{h}[S(p_2,p_3,x)+S(p_1,p,x)-px]}\frac{d^dpd^dx}{(2\pi h)^{d/2}} $$
be equal. This allows us to provide in Section \ref{kontform} a quick but non rigorous proof of the existence of
the generating function relying only on the associativity of the Kontsevich star product.

\par

The second way to derive the SGA equation is symplectic reduction. Consider the symplectic groupoid
$G$ over $M=\mathbb R^d$ as above. Let us call $\mathcal L_S\subset \overline G\times \overline G\times G$ 
the Lagrangian submanifold associated to the generating function $S$ (i.e., $\mathcal L_S = \graph (dS)$). 
Now consider the spaces $H(k) = \overline G^k\times G$ and the diagonal $\Delta_{l_1,\dots,l_k}\subset
\overline{H(k)}\times H(l_1)\times \dots\times H(l_k)$, 
$$\Delta_{l_1,\dots,l_k} = \big\{ (g_1,\dots,g_k,y),(x_{11},\dots,x_{1l_1},g_1),\dots,(x_{k1},\dots,x_{kl_k},g_k)\big\}.$$
This is a coisotropic subspace of $H(k)\times H(l_1)\times\dots\times H(l_k)$. Then one can consider the
symplectic reduction by the diagonal $\Delta_{l_1,\dots,l_k}$ which sends Lagrangian submanifolds of 
$H(k)\times H(l_1)\times\dots\times H(l_k)$ to Lagrangian submanifolds of $H(l_1+\dots+l_k)$. In particular
$\mathcal L_S\oplus\mathcal L_S\oplus \mathcal L_I$ ($I(p,x)= px$) is sent to $\mathcal L_1 \subset H(3)$ and
$\mathcal L_S\oplus\mathcal L_I\oplus \mathcal L_S$ to $\mathcal L_2\subset H(3)$. One can check that $\mathcal L_1=\mathcal  L_2$
iff $S$ satisfies the SGA equation. In fact we have here, hidden in the background, a structure of an operad,
 the Lagrangian operad (see \cite{Cattaneo2002}).
 
%
%
%
%
%

Now consider $M= \mathbb R^d$ with the zero Poisson structure.  The symplectic groupoid $G_0$ over it
is the cotangent bundle ($G_0 = \mathbb{R}^{*d}\times\mathbb{R}^d$). The source map  and
the target map $s,t:G_0\rightarrow \mathbb R^d$ are identified with the cotangent bundle projection.
The inclusion $\epsilon:\mathbb R^d \rightarrow G_0$ is defined by $\epsilon(x) = (0,x)$, the inverse
map $i:G_0\rightarrow G_0$ by $i(p,x)=(-p,x)$ and the product is the fiber wise addition, i.e.,  $(p_1,x)\bullet(p_2,x) =
(p_1+p_2,x)$.  The product space $G_0^{(m)}$  can be seen
as the graph of the differential of the function $S_0(p_1,p_2,x)= x(p_1+p_2)$. It is
easy to check that $S_0$ satisfies the SGA equation. We investigate deformations of this
trivial generating function. Let us be more precise.

\par
\begin{Definition}
A \textbf{deformation} of the trivial generating function is a formal
power series in $h$,  $S_h = S_0 +hS_1+h^2S_2\dots$, obeying the SGA equation and
such that $S_0(p_1,p_2,x)=x(p_1+p_2)$.
\par
Such a deformation  is called \textbf{natural} if 
\begin{enumerate}
\item $S_n(p,q,x)$ are polynomial in $q$ and $q$
\item $S_n(\lambda q,\lambda p,x) = \lambda^{n+1} S_n(p,q,x)$
\item $S_n(p,0,x)=S_n(0,p,x)=0$
\item $S_n^i(p,p) = 0$, where $S_n^i$ is the homogeneous part  of $S_n$ of degree $i$ in the first argument.
\end{enumerate}

\end{Definition}

In Section \ref{recovering} we show that, provided we 
have a  natural deformation $S_h=S_0 +hS_1+h^2S_2+\dots$ of the trivial generating function,
we can deform the structure maps of the trivial symplectic groupoid into
\[
\begin{array}{cccc}
\epsilon_h(x) & = & (0,x)  &\textbf{unit map}\\
i_h(p,x)      & = & (-p,x) &\textbf{inverse map}\\
s_h(p,x)      & = & \nabla_{p_2}S_h(p,0,x) & \textbf{source map}\\
t_h(p,x)      & = & \nabla_{p_1}S_h(0,p,x) & \textbf{target map}
\end{array}
\]
such that the groupoid structure is (formally) preserved. 
\par
Moreover  there is a unique Poisson bracket on $\mathbb R^d$  such that the source, $s_h$, is a Poisson
map with respect to the canonical symplectic structure on the formal symplectic groupoid. 
This Poisson bracket is given by $\{f,g\}_{\mathbb R^d}(x)=2hS_1(df,dg,x)$, the first order term of
the generating function.
We can now formulate the formal integration problem for Poisson manifolds.

\vspace{1cm}
\textbf{Formal integration problem for Poisson manifolds:}
\par

\emph{
Given a Poisson structure on $\mathbb R^d$, does there exists a deformation of the trivial generating
function such that the  first order term  is the original Poisson structure?
}

\subsection{Main Result, main example, main interpretation}

The following theorem gives a positive answer to the deformation problem for symplectic
groupoids. This is the main result of this article.

\begin{Theorem}\label{th:main}
Given a Poisson structure $\alpha$ on $\mathbb R^d$ 
there exists a unique natural deformation of the trivial generating function such that
the first order is precisely $\alpha$. Moreover we have an explicit formula for this deformation 
$$S_h(p_1,p_2,x) = x(p_1+p_2)+\sum_{n=1}^\infty \frac{h^n}{n!}\sum_{\Gamma\in T_{n,2}}W_\Gamma \hat B_\Gamma(p_1,p_2,x)$$
where $T_{n,2}$ is the set of Kontsevich trees of type $(n,2)$, $W_\Gamma$ is the
Kontsevich weight of $\Gamma$ and $\hat B_\Gamma$ is the symbol of the bidifferential
operator $B_\Gamma$ associated to $\Gamma$.
\end{Theorem}


Section \ref{recovering} explains how to recover the structure maps from the deformed generating function.
In Section \ref{examples} we  present basic examples of formal symplectic groupoids.
In particular the main one is in the case of a linear Poisson structure $\alpha^{ij}(x) = \alpha^{ij}_k x^k$, i.e., 
when one considers the Kirillov--Kostant Poisson structure on the dual $\mathcal G^*$ of a Lie algebra $\mathcal G$.
In this case, the generating function of the symplectic groupoid over 
$\mathcal G^*$ reduces exactly to the Campbell--Baker--Hausdorff formula$$
S_h(p_1,p_2,x) = \langle \frac{1}{h}CBH(h p_1,h p_2), x\rangle
$$
where $\langle,\rangle$ is the natural pairing between $\mathcal G$ and $\mathcal G^*$.

This basic example suggests to consider the generating function as a  generalized CBH formula to the non-linear
case and reproves in the linear case a result of V. Kathotia(\cite{kathotia1998}).

\par

Sections \ref{Cayleytree}  to \ref{proof} are  devoted to the proof of the  Theorem \ref{th:main}. In Section \ref{Cayleytree} we 
introduce special graphs, the Cayley trees, which allow us to write down a perturbative version of the SGA equation. 
In Section \ref{sec:Konttrees}  we describe the Kontsevich trees. We use them to produce an
explicit solution for the deformation problem. Section \ref{proof} completes the proof of Theorem \ref{th:main}.

In the last section  we come to the comparison with deformation quantization.
We see that the Kontsevich star-product can be put into the form
$$
f*g(x) = \exp\Big(\frac{1}{h}\sum_{i=0}^{\infty}h^i D_i(h\partial_y,h\partial_z,x)\Big)f(y)g(z)_{\Big|y=z=x}
$$
where $$ S_h(p_1,p_2,x) =x(p_1+p_2)+\frac{1}{h} D_0(h p_1,h p_2,x).$$ 
This allows us to interpret the generating function as a semi-classical version of the Kontsevich star-product formula.
At last, considering associativity of star product of exponential functions, we are able to provide an elegant but
non-rigourous proof of the existence part of Theorem \ref{th:main}.

\subsection{Planned developements}

One of the next objective is to carry the construction of the formal symplectic
groupoid to a general Poisson manifold. Such a globalisation has already been made
by Karabegov \cite{karabegov2003} in the symplectic case.

A second possible developement is to try to derive the existence of the 
deformation of the trivial generating function from a kind of ``semi-classical"
formality theorem. 

At last we plan to compare the formal construction carried out in this article with the 
non-formal construction coming from the Poisson-sigma model (see \cite{feldcat2000}) and
with the local symplectic groupoid construction of \cite{CDW1987} and \cite{karasev1987}.

\begin{Ack}
The second author thanks Ernst Hairer for useful discussions, and suggestions.
\end{Ack}

\section{Recovering the formal groupoid from the generating function} \label{recovering}

In this section we show that one can recover formally the structure of 
symplectic groupoid from a generating function obeying the SGA equation.

\begin{Proposition}\label{formalgroupoid}
Let $S_h$ be a natural deformation of the trivial generating function which
satisfies the SGA equation. Then the set $G_h = \mathbb R^{*d}[[h]]\times \mathbb R^d[[h]]$
can be given a structure of \textbf{formal symplectic groupoid}, i.e., the maps
\[
\begin{array}{cccc}
\epsilon_h(x) & = & (0,x)  &\textbf{unit map}\\
i_h(p,x)      & = & (-p,x) &\textbf{inverse map}\\
s_h(p,x)      & = & \nabla_{p_2}S_h(p,0,x) & \textbf{source map}\\
t_h(p,x)      & = & \nabla_{p_1}S_h(0,p,x) & \textbf{target map}
\end{array}
\]
and the multiplication given by 
$$G_h^{(m)}= \graph(dS)$$
satisfy formally the axioms of a groupoid.
\par
In particular, if we endow $G_h$ with the canonical symplectic form, then 
$G_h^{(m)}$ is formally Lagrangian in $\overline{G_h}\times \overline{G_h}\times G_h$.
\end{Proposition}

\begin{proof}
The multiplication space being given by the graph of the differential
of the generating function, we have automatically that the product, when
defined, is associative (it satisfies the SGA equation) and that $G_h^{(m)}$ is
formally a Lagrangian submanifold of $\overline{G_h}\times\overline{G_h}\times G$. We still
have to check that the space of composable pairs is the right one, i.e., $(g,h)\in G_h^{(2)}$ iff
$s(g) = t(h)$. We do that by noticing that all products are of the form
$(p_1,\nabla_{p_1} S_h(p_1,p_2,x))\bullet (p_2,\nabla_{p_2} S_h(p_1,p_2,x)) = (\nabla_x S_h(p_1,p_2,x),x)$.
Thus the check amounts to see that
$$s(p_1,\nabla_{p_1}S_h(p_1,p_2,x))=t(p_2,\nabla_{p_2}S_h(p_1,p_2,x))$$
which can be seen by differentiating the SGA equation with respect to $p_2$, putting $p_2 = 0$ and using the fact
that $S_h$ is natural.

It remains still  to check the following axioms
\begin{gather*}
t(gh)=t(g)\quad (1),
\quad 
s(gh)=s(h)\quad (2),
\quad
\epsilon(t(g))g=g\quad (3),
\quad g\epsilon(s(g))=g\quad (4),
\\
s(i(g))=t(g)\quad (5),
\quad
i(g)g=\epsilon(s(g))\quad (6),
\quad
gi(g)=\epsilon(t(g))\quad (7),
\end{gather*}

Axiom 1 is obtained by differentiating the SGA equation w.r.t. $p_1$, putting $p_1=0$ and using
naturality of $S_h$. The Axiom 2 is similar but for replacing $p_1$ by $p_3$. Axiom 3 and Axiom 4
are direct consequences of the naturality. The last three axioms are however a bit more
tricky. First let us prove two Lemmas.

\begin{Lemma}[\bf{Inversion of source and target}]\label{formaldiffeo}
Denote $F_p(x) = \nabla_{p_2}S_h(p,0,x)$ and $G_p(x) = \nabla_{p_1}S_h(0,p,x)$. Then $F_p$
and $G_p$ are formal diffeomorphisms and their inverses are given by 
$$F^{-1}_p(x) = \nabla_{p_2}S_h(-p,p,x),\quad G^{-1}_p(x) = \nabla_{p_1}S_h(p,-p,x).$$
\end{Lemma}

\begin{proof}
Denote $\overline F_p(x) = \nabla_{p_2}S_h(-p,p,x)$ and $\overline G_p(x) = \nabla_{p_1}S_h(p,-p,x)$.
Differentiating the SGA equation w.r.t. $p_1$, putting $p_1=p$, $p_2=-p$, $p_3=p$, we get that $\overline G_p\circ G_p = id$.
Putting $p_1=0$, $p_2=p$, $p_3=-p$, we get $G_p\circ \overline G_p = id$. Thus $\overline G_p = G^{-1}_p$.
Similarly differentiating the SGA equation w.r.t. $p_3$, putting $p_1=p$, $p_2=-p$, $p_3=p$, we get that $\overline F_p\circ F_p = id$.
Putting $p_1=p$, $p_2=-p$, $p_3=0$, we get $F_p\circ \overline F_p = id$. Thus $\overline F_p = F^{-1}_p$.
\end{proof}

\begin{Lemma}[\bf{Relation between source and target}]\label{relationformdiffeo}
Denote $F_p(x) = \nabla_{p_2}S_h(p,0,x)$ and $G_p(x) = \nabla_{p_1}S_h(0,p,x)$. Then we
have the relation $$F_p = G_{-p}.$$
\end{Lemma}

\begin{proof}
Notice that it is equivalent to prove that $F_p = G_{-p}$ or $F^{-1}_p = G^{-1}_{-p}$.
We prove the second identity. For each $n\geq 1$ we have the decomposition
$$S_n(p_1,p_2,x) = S_n^1(p_1,p_2,x)+  S_n^2(p_1,p_2,x)+\dots+ S_n^n(p_1,p_2,x)$$
where $S_n^i$ is the part of $S_n$ which is homogeneous of degree $i$ in the first argument.
Now we have that $$S_n^i(-p,p,x)= (-1)^iS_n^i(p,p,x) = 0$$ because of naturality of 
the generating function. This implies that $S(-p,p,x) = 0$. If we differentiate this equation
with respect to p we get exactly $F^{-1}_p = G^{-1}_{-p}$.

\end{proof}

Going back to the check of axioms we get that Axiom 5 is exactly equivalent to $F_p = G_{-p}$.
\par
As for Axiom 6, if we pose $i(p,x) = (-p,\nabla_{p_1}S_h(-p,p,s(p,x))$, then 
\begin{eqnarray*}
i(p,x)(p,x) & = & (\nabla_xS_h(-p,p,s(p,x),s(p,x))\\
            & = & (0,s(p,x))\\
            & = &\epsilon(s(p,x)),
\end{eqnarray*}
provided that $x = \nabla_{p_2}S_h(-p,p,\nabla_{p_2}S_h(p,0,x))$, which is guaranteed by Lemma \ref{formaldiffeo}.

Similarly for Axiom 7, if we put $\tilde i(p,x)= (-p,\nabla_{p_2}(p,-p,t(p,x))$ we get that $(p,x)\tilde i(p,x) = \epsilon(t(p,x))$.
Now by Lemma \ref{relationformdiffeo} we get
$$ \nabla_{p_2}S_h(p,-p,t(p,x)) = F_{-p}\circ G_p(x) = x$$
$$ \nabla_{p_1}S_h(-p,p,s(p,x)) = G_{-p}\circ F_p(x) = x. $$
Thus $i(p,x) = \tilde i(p,x) = (-p,x)$.
\end{proof}

Now using the canonical symplectic bracket on $G_h$, i.e., 
$$\{F,G\}_{G_h}(p,x)= \langle\nabla_x F(p,x),\nabla_p G(p,x)\rangle-\langle\nabla_p F(p,x),\nabla_x G(p,x)\rangle$$
we can consider the problem of finding a Poisson bracket on $\mathbb R^d$ such that $s_h$ is Poisson algebra homomorphism, i.e.,
$$s^*_h\{f,g\}_{\mathbb R^d}(p,x) = \{s_h^*f,s_h^*g\}_{G_h}(p,x).$$
The following Proposition answers this question.

\begin{Proposition}
There is a unique Poisson structure on $\mathbb R^d$ such that $s^*$ is a (formal) Poisson map. Moreover this
Poisson structure is given by $$\{f,g\}_{\mathbb R^d}(x) = \{s_h^*f,s_h^*g\}_{G_h}(0,x)= 2hS_1(df,dg,x).$$
\end{Proposition}

\begin{proof}
Suppose there exists a Poisson structure $\{,\}_{\mathbb R^d}$ such that $s^*$ is Poisson. 
This means that $\{f,g\}_{\mathbb R^d}(s_h(p,x)) = \{s_h^*f,s_h^*g\}_{G_h}(p,x)$. In particular
if we put $p=0$ we get exactly that $\{f,g\}_{\mathbb R^d}(x)  = \{s_h^*f,s_h^*g\}_{G_h}(0,x)= 2hS_1(df,dg,x)$ which
show uniqueness.
Now it remains to prove that $\{f,g\}_{\mathbb R^d}(s_h(p,x)) = \{s_h^*f,s_h^*g\}_{G_h}(p,x)$ which proves 
as well that the induced bracket is Poisson.

\par
Then we have to check that $$\{s_h^*f,s_h^*g\}_{G_h}(0,s_h(p,x)) = \{s_h^*f,s_h^*g\}_{G_h}(p,x).$$
An easy computation gives us that this equation is equivalent to the following
$$\frac{\partial s_h^k}{\partial p^l}(0,s_h(p,x))-\frac{\partial s_h^l}{\partial p^k}(0,s_h(p,x)) =
\sum_{i=1}^d\big\{\frac{\partial s_h^k}{\partial x^i}(p,x) \frac{\partial s_h^l}{\partial p^i}(p,x)-
\frac{\partial s_h^l}{\partial x^i}(p,x) \frac{\partial s_h^k}{\partial p^i}(p,x)\big\}. $$
Differentiating the SGA equation first with respect to $p_3$ and then to $p_2$ and then putting $p_1=p$, $p_2=p_3=0$,
we get
$$\nabla_{p^1_k}\nabla_{p^2_l}S(0,0,s_h(p,x)) = \sum_{i=1}^d\nabla_{x^i}\nabla_{p^2_k}S(p,0,x)\nabla_{p^1_i}
\nabla_{p^2_l}S(p,0,x)-\nabla_{p^2_k}\nabla_{p^2_l}S(p,0,x).$$
Taking the difference between this equation and the same but with the indices $k$ and $l$ interchanged we
finish the proof.
\end{proof}

\section{Basic examples}\label{examples}

Let us see on some examples what are the generating functions and the formal symplectic groupoids.
We already know what happens in the case of the trivial Poisson structure over $\mathbb R^d$. The 
generating function is $S_0(p_1,p_2,x) = (p_1+p_2)x$ and the associated symplectic groupoid is
the cotangent bundle $T^*\mathbb R^d$ with structure maps $s(p,x) = x$, $t(p,x)=x$, $\epsilon(x) = (0,x)$
$i(p,x) = (-p,x)$. The composition is the fiberwise addition.

\subsubsection{Constant Poisson structure}
Suppose one has a constant Poisson structure $\alpha(x) = \alpha$. 
The main Theorem tells us that the generating function is $S_h(p_1,p_2,x)= S_0(p_1,p_2,x) + h p_1^t\alpha p_2$. 
The multiplication space can then be described as 
$$ G^{(m)}_h=\Big\{\Big((p_1,x+h\alpha p_2),(p_2,x-h\alpha p_1),(p_1+p_2,x)\Big),(p_1,p_2,x)\in B_2\Big\}.$$
By Proposition \ref{formalgroupoid} the structure maps are given by

\begin{eqnarray*}
\epsilon_h(x) & = & (0,x)\\
i_h(p,x) & = & (-p,x)\\
s_h(p,x) &=& x-h \alpha p\\
t_h(p,x)& =& x+h\alpha p.
\end{eqnarray*}

\subsubsection{Linear Poisson structure}

Suppose that we have a linear Poisson structure $\alpha^{ij}(x)=\alpha^{ij}_kx^k$ on $\mathbb R^d$ which
can then be considered as the dual of a Lie algebra $\mathbb R^d = \mathcal G^*$, the bracket on $\mathcal G$
being given by $[\epsilon^i,\epsilon^j] = 2\alpha^{ij}_k \epsilon^k$, where $\epsilon^l$, $l=1,\dots,d$ is a basis of
$\mathcal G (= \mathbb R^{*d})$. For this Lie algebra we have the CBH formula $\exp(p_1)\exp(p_2)=\exp(CBH(p_1,p_2))$,
$$CBH(p_1,p_2)= p_1+p_2+\frac{1}{2}[p_1,p_2]
+\frac{1}{12}\Big([p_1,[p_2,p_2]]+[p_2,[p_2,p_1]]\Big)+\dots,$$

It is easy to check directly that 
$$ S_h(p_1,p_2,x) = \langle\frac{1}{h}CBH(h p_1,h p_2),x\rangle $$
where $\langle.,.\rangle$ is the usual pairing between $\mathcal G$ and $ \mathcal G^*$,
satisfies the SGA equation. It is equivalent to the associativity of CBH, i.e.,
$$
CBH(hp_1,CBH(hp_2,hp_3)) = CBH(CBH(hp_1,hp_2),hp_3).
$$

By the uniqueness of the generating function given by the main Theorem we recover 
a result of V. Kathotia (see \cite{kathotia1998}):

\begin{Proposition} For the Poissson structure coming from the dual of a Lie algebra we have
$$\langle\frac{1}{h}CBH(h p_1,h p_2)-(p_1+p_2),x\rangle = \sum_{n\geq1}\frac{h^n}{n!}\sum_{\Gamma\in T_{n,2}}W_\Gamma \hat B_\Gamma(p_1,p_2,x)$$ 
\end{Proposition}

This result is one of the main ingredient
to prove that CBH-quantization is a deformation quantization in the case
of the dual of a Lie algebra. It allows us to consider the generating function
as a generalization of the CBH formula to the non-linear case.

By Proposition \ref{formalgroupoid} we have that the deformed source and target maps are
$$t_h(p,x) = \langle\frac{1}{h} \nabla_{p_1} CBH(0,h p),x\rangle = x +h \alpha^{ij}_k x^k p_j+\frac{h^2}{3} 
\alpha^{uv}_l\alpha^{ni}_v x^l p_u p_n+\dots $$
$$s_h(p,x) = \langle\frac{1}{h} \nabla_{p_2} CBH(h p,0),x\rangle = x -h \alpha^{ij}_k x^k p_j+\frac{h^2}{3} 
\alpha^{uv}_l\alpha^{ni}_v x^l p_u p_n+\dots $$

\section{Perturbative form of the SGA equation} \label{Cayley}

\label{Cayleytree}

The goal of this section is to formulate a perturbative version of the SGA equation. It is
divided in two parts. First we introduce some tools and state the perturbative version of the SGA equation in
Proposition \ref{pertSGA}. The proof is then split into several lemmas. 

\subsection{Perturbative SGA and Cayley trees}

Let us recall that $B_n:= (\mathbb R^{*d})^n \times \mathbb R^d$. First suppose that we are looking for
a generating function of the form $S_h = S_0+h S$ where $S_0(p_1,p_2,x) = (p_1+p_2)x$ is 
the trivial generating function and $S\in C^\infty (B_2)[[h]] $ is a formal series $S = S_1 +h S_2 +\dots$.
Inserting $S_h$ in the SGA equation we get a new version of this equation for $S$, $M^1(S) = M^2(S)$ 
where $M^i:C^\infty (B_2)[[h]]\rightarrow C^\infty (B_3)[[h]]$ are defined
by 

\begin{eqnarray*}
M^1(S)(p_1,p_2,p_3,x;h)  & = & h S(p_{1},p_{2},\bar{x})+h S(\bar{p},p_{3},x) \\
                             &   &  \quad\quad- h^2 \nabla_xS(p_1,p_2,\bar{x})\nabla_{p_1}S(\bar{p},p_3,x)
\end{eqnarray*}
\begin{eqnarray*}
 \bar{p} & = &p_1+p_2+h\nabla_x S(p_{1},p_{2},\bar{x}),\\
\quad \bar{x} & = &x+h\nabla_{p_1}S(\bar{p},p_{3},x) 
\end{eqnarray*}
and
\begin{eqnarray*}
M^2(S)(p_1,p_2,p_3,x;h)  & = & h S(p_{2},p_{3},\tilde{x})-\tilde{p}\tilde{x} +h S(p_{1},\tilde{p},x)\\
                     &   &  \quad\quad- h^2 \nabla_xS(p_2,p_3,\tilde{x})\nabla_{p_2}S(p_1,\tilde{p},x)
\end{eqnarray*}
\begin{eqnarray*}
\tilde{p} & = &p_2+p_3+h\nabla_x S(p_{2},p_{3},\tilde{x}),\\
\quad \tilde{x} & = &x+h\nabla_{p_2}S(p_{1},\tilde{p},x).
\end{eqnarray*}

\par

The idea now is to expand $M^i(S)(h)$, $i=1,2$ into powers of $h$ and then to analyze the conditions
imposed on $S$ by the equation at each order. For that purpose we are going to introduce some tools and
methods that are heavily inspired from the tools used in numerical analysis to determine the order
condition of a Runge--Kutta method. The main ingredients are trees which are used to represent
the so called elementary differentials and elementary functions. As these ideas go back to Cayley, we
 call such trees Cayley trees, in order to distinguish them from Kontsevich trees which will
also appear in the story. In the sequel we will mainly  follow the notations of \cite{GeomInt}.  
\begin{Definition}
        \begin{enumerate}
                \item{A \textbf{graph} $t$ is given by a set of vertices $V_t=\{1,\dots,n\}$ and
                        a set of edges $E_t$ which is a set of pairs of elements of $V_t$. We denote the number
                        of vertices by $|t|$. An \textbf{isomorphism} between two graphs $t$ and $t'$ having
                        the same number of vertices is 
                        a permutation $\sigma\in S_{|t|}$ such that $\{\sigma(v),\sigma(w)\}\in E_{t'}$ if
                        $\{v,w\}\in E_{t}$. Two graphs are called \textbf{equivalent} if
                        there is an isomorphism between them. The \textbf{symmetries} of a graph are
                        the automorphisms of the graph. We denote the group of symmetries by $sym(t)$.}

                \item{A \textbf{tree} is a graph which has no cycles. Isomorphisms and symmetries are
                        defined the same way as for graphs }

                \item{A \textbf{rooted tree} is a tree  with one distinguished vertex. An 
                        \textbf{isomorphism} of rooted trees is an isomorphism of graphs which sends the
                        root to the root. Symmetries and equivalence are defined correspondingly.}

                \item{A \textbf{bipartite graph} is a graph $t$ together with a map $\omega:V_t\rightarrow
                        \{\circ,\bullet\}$ such that $\omega(v)\neq\omega(w)$ if $\{v,w\}\in E_t$. An isomorphism
                        of bipartite trees is an isomorphism of graphs which respects the coloring, i.e., $\omega(\sigma
                        (v))=\omega(v)$.}
        \end{enumerate}
\end{Definition}

The following table summarizes some notations we will use in the sequel.

\begin{tabular}{|c|l|}
\hline
$T$             & the set of bipartite trees\\
$RT$            & the set of rooted bipartite trees\\
$RT_\circ$      & the set of elements of $RT$ with white root\\
$RT_\bullet$    & the set of elements of $RT$ with black root\\
\hline
\end{tabular}

$[A]$       : the  set of equivalence classes of graphs in $A$ (ex: $[RT]$). They are called
\textbf{topological} ``$A$" trees.
\par

The elements of $[RT]$ can be described recursively as follows
\par
\begin{enumerate}
        \item{$\circ,\bullet\in [RT]$}
        \item{if $t_1,\dots,t_m\in [RT_\circ]$, then the tree $[t_1,\dots,t_m]_\bullet\in [RT]$ where
                 $[t_1,\dots,t_m]_\bullet$ is  defined by 
                connecting the roots of $t_1,\dots,t_m$ with $\bullet$ and saying that $\bullet$ is
                the new root. And the same if we interchange $\circ$ and $\bullet$. }
\end{enumerate}

Now with the help of this recursive description of topological rooted trees we define   elementary
differentials and elementary generating functions.

\begin{Definition}[\textbf{Elementary Differentials (ED)}]
Let $i=1,2$, $t\in[RT]$. The elementary differential $D^iS(t)$ of $S\in C^\infty(B_2)[[h]]$ is defined recursively
as follows,

\begin{enumerate}
\item $D^iS(\circ) = \nabla_x S$ , $D^iS(\bullet) = \nabla_{p^i} S$
\item $D^iS(t) = \nabla_{p^i}^{(m+1)}S(D^iS(t_1),\dots,D^iS(t_m))$ if $t= [t_1,\dots,t_m]_\bullet$
\item $D^iS(t) = \nabla_x^{(m+1)}S(D^iS(t_1),\dots,D^iS(t_m))$ if $t= [t_1,\dots,t_m]_\circ$
\end{enumerate}
where $\nabla_x^{(k)}S $ stands for the $k^{th}$ derivative of $S$ w.r.t. $x$ evaluated at $(p_1,p_2,x)$ if $i=1$ and at
$(p_1+p_2,p_3,x)$ if $i=2$.  $\nabla_{p^i}^{(k)}S$ stands for the $k^{th}$ derivative of $S$ w.r.t. $p^i$ evaluated at  $(p_1+p_2,p_3,x)$ 
if $i=1$ and at $(p_1,p_2+p_3)$ if $i=2$.
\end{Definition}

\begin{Definition}[\textbf{Elementary Generating Functions (EGF)}]
Let $i=1,2$, $t\in[RT]$. The elementary generating function $S^i(t)$ of $S\in C^\infty(B_2)[[h]]$ is defined recursively
as follows,
\begin{enumerate}
\item $S^1(\circ) =  S(p_1,p_2,x)$ ,\quad $S^1(\bullet) =  S(p_1+p_2,p_3,x)$
\item $S^2(\circ) =  S(p_2,p_3,x)$ ,\quad $S^2(\bullet) =  S(p_1,p_2+p_3,x)$
\item $S^i(t) = \nabla_{p^i}^{(m+1)}S(D^iS(t_1),\dots,D^iS(t_m))$ if $t= [t_1,\dots,t_m]_\bullet$
\item $S^i(t) = \nabla_x^{(m+1)}S(D^iS(t_1),\dots,D^iS(t_m))$ if $t= [t_1,\dots,t_m]_\circ$
\end{enumerate}
with the same notation as above.
\end{Definition}

Some examples are given in the following table:

\vspace{0.5cm}
\begin{tabular}{|c|c|c|c|}
\hline
Diagram & Notation            & ED                           &  EGF \\
&&&\\
\hline
\rtree{RT-1}&$[\bullet]_\circ $ & $\nabla^{(2)}_x S\nabla_p S $ & $\nabla_x S\nabla_p S$ \\ 
&&&\\
\rtree{RT-2}&$[\circ,\circ]_\bullet$ &$ \nabla^{(3)}_pS(\nabla_x S, \nabla_x S) $& $\nabla^{(2)}_pS(\nabla_x S, \nabla_x S)$  \\ 
&&&\\
&&&\\
\rtree{RT-3}&$[\bullet,[\circ]_\bullet]_\circ$ &$ \nabla^{(3)}_xS(\nabla_p S,\nabla^{(2)}_p S\nabla_x S ) $ & $ \nabla^{(2)}_xS(\nabla_p S,\nabla^{(2)}_p S\nabla_x S ) $   \\
&&&\\
\hline
\end{tabular}
\vspace{0.5cm}

Remark that for EGF it is not important which vertex is the root. This is not the case for ED. Let us be more
precise.
\begin{Definition}[\textbf{Butcher product}]
Let $u = [u_1,\dots,u_k], v = [v_1,\dots,v_l]\in [RT]$.
We denote by
\begin{eqnarray*}
u\circ v &=& [u_1,\dots,u_k,v]\\
v\circ u &=& [v_1,\dots,v_l,u]\\
\end{eqnarray*}
the Butcher product. We have not written the obvious conditions
on the $u_i$ and $v_i$ so that the product remains bipartite.
\end{Definition}

\begin{Definition}[\textbf{Equivalence relation on rooted topological trees}]
We consider the minimal equivalence relation on $[RT]$ such that
$u\circ v \sim  v\circ u$. 
\end{Definition}

\textbf{Properties of this relation:}
\par
It is clear that
\begin{enumerate}
\item two topological rooted trees are equivalent if it is possible to pass from one to the
other by changing the root. More precisely:  $t,t'\in[RT]$, $t\sim t'$ iff there exists a 
representative $(E,V,r)$ of $t$ and a representative $(E',V',r')$ of $t'$ and a vertex $r''\in V$ such
that $(E,V,r'')$ and $(E',V',r')$ are isomorphic rooted trees.
\item the quotient of $[RT]$ by this equivalence relation is exactly $[T]$.
\item it follows immediately from the definition $S^i(t) = S^i(t')$ if $t\sim t'$ for $i=1,2$.
\end{enumerate}

Then, it makes sense to define the EGF on bipartite trees.

\begin{Definition}
Let $S\in C^\infty (B_2)[[h]]$ and  $t=(V_t,E_t)\in T$. Then

$$S^1(t) := \sum_{\beta:E_t\rightarrow \{ 1,\dots,d\}}\prod_{v\in V_t}\Big[\Big(
\prod_{\substack{e\in E_t\\e=\{*,v\}}}D^{1,\omega(v)}_{\beta(e)} \Big)S \Big]$$
where
$$D^{1,\bullet}_{\beta(e_1)}\dots D^{1,\bullet}_{\beta(e_k)}S := \frac{\partial^k S}{\partial p^{1}_{\beta(e_1)}\dots 
\partial p^{1}_{\beta(e_k)}}(p_1+p_2,p_3,x)$$

$$D^{1,\circ}_{\beta(e_1)}\dots D^{1,\circ}_{\beta(e_k)}S := \frac{\partial^k S}{\partial x^{\beta(e_1)}\dots \partial x^{\beta(e_k)}}(p_1,p_2+p_3,x)$$
and correspondingly for $S^2(t)$.
 \end{Definition}

It is clear that this new definition of $S^i(t)$ is equivalent to the previously introduced recursive one.
This definition is however better if we want to deal with the fact that $S$ is a formal series. Namely we immediately get 
the relation

\begin{eqnarray*}
h^{|t|}S^i(t) & = & h^{|t|}\sum_{\beta:E_t\rightarrow\{1,\dots,d\}}\prod_{v\in V_t}\Big[\big(%
          \prod_{\substack{e\in E_t\\e=\{*,v\}}}D^{i,\omega(v)}_{\beta(e)}\big) \big\{ \frac{1}{h}%
          \sum_{n=1}^{\infty}h^n S_n\big\}\Big]\\
          & = & \sum_{n=|t|}^{\infty}h^n\sum_{\substack{n_1+\dots+n_{|t|}=n\\n_i\geq 1}}C^i_t(S_{n_1},\dots,S_{n_{|t|}})
\end{eqnarray*}

which defines the $C^i_t$ which are multi-differential maps from $C^\infty(B_2)^{|t|}$ to $C^\infty(B_3)$.
\par
We can now state the main Proposition of this section.

\begin{Proposition}[\textbf{Perturbative version of the SGA equation}]\label{pertSGA}
The formal series $S_h= S_0 +\sum_{n\geq 1}h^n S_n$ satisfies the SGA equation iff
for each $n>0$ we have 

$$
\sum_{\substack{t\in T\\ |t|\leq n}}\frac{1}{|t|!}\sum_{\substack{n_1+\dots+n_{|t|}=n\\n_i\geq 1}}C_t^1(S_{n_1},\dots,S_{n_{|t|}})
-C_t^2(S_{n_1},\dots,S_{n_{|t|}})=0.
$$
\end{Proposition}

Let us remark that for all $f\in C^\infty(B_2)$ we have that,
$$C^1_\bullet(f)+C^1_\circ(f)-C^2_\bullet(f)-C^2_\circ(f) = dS$$
where $d:C^\infty (B_n)\rightarrow C^\infty(B_{n+1})$ is a differential 
(i.e., $d^2=0$) defined by the formula

$$d f(p_{1},\dots,p_{n+1}) =
f(p_{2},...,p_{n+1})-\sum_{i=1}^{n}(-1)^{i+1}f(p_{1},\dots,p_{i}+p_{i+1})
$$$$+(-1)^ {n+1}f(p_{1},...,p_{n}).$$
This differential can be interpreted either as the Hochschild differential on symbols
of multi-differential operators on $C^\infty(\mathbb R^d)$ or as the differential of
the trivial symplectic groupoid cohomology over $\mathbb R^d$. This remark allows us
to put the previous recursive equations into the form $$ dS_n + H_n(S_{n-1},\dots,S_1)=0$$
which is exactly the analog of the recursive equation involved when considering
star-products.
\par
The remaining of this section is devoted to proving Proposition \ref{pertSGA}.

\subsection{Proof of the Proposition}\label{defeqn}

It follows from a series of little Lemmas.
We are first interested in expanding 
 \begin{eqnarray}
 \bar{p} &=&p_1+p_2+h\nabla_x S(p_{1},p_{2},\bar{x}),\\
 \bar{x} &=&x+h\nabla_{p_1}S(\bar{p},p_{3},x), 
 \end{eqnarray}
and 
 \begin{eqnarray}
\tilde{p} &=&p_2+p_3+h\nabla_x S(p_{2},p_{3},\tilde{x}),\\
\tilde{x} &=&x+h\nabla_{p_2}S(p_{1},\tilde{p},x).
 \end{eqnarray}
as power series in $h$.

The method used is essentially the same as in numerical analysis when one wants
to express the Taylor series of the numerical flow of  a Runge--Kutta method.
Namely the equations above have a form very close to the partioned implicit Euler method(see \cite{GeomInt}).

\begin{Definition}
Let $t = [t_1,\dots,t_m]\in [RT]$. Consider the list $\tilde t_1,\dots,\tilde t_k$ of all non isomorphic
trees appearing in $t_1,\dots,t_m$. Define $\mu_i$ as the number of time the tree $\tilde t_i$ appears in
$t_1,\dots,t_m$. Then we introduce the \textbf{symmetry coefficient} $\sigma(t)$  of 
$t$ by the following recursive definition:
\par
$$\sigma(t) = \mu_1!\mu_2!\dots\sigma(\tilde t_1)\dots\sigma(\tilde t_k).$$
Moreover $\sigma(\circ)=\sigma(\bullet)= 1$.
\end{Definition}

It is clear that $\sigma(t)$ is the number of symmetries for each representative of $t$ (i.e $\sigma(t) = |Sym(t')|$ for
all $t'\in t$).

\begin{Lemma}

 There exist unique formal series for $\bar x,\bar p$ (resp. $\tilde x$, $\tilde p$) which
 satisfy equation (1) and (2) (resp. (2) and (3)). They are given by
 \begin{eqnarray}
 \bar{x}(h) &=&x+\sum_{t\in [RT_{\bullet}]}\frac{h^ {|t|}}{\sigma(t)}D^1S(t),\\
 \bar{p}(h) &=&p_1+p_2+\sum_{t\in [RT_{\circ}]}\frac{h^ {|t|}}{\sigma(t)}D^1S(t),
 \end{eqnarray}
and by
 \begin{eqnarray}
 \tilde{x}(h) &=&x+\sum_{t\in [RT_{\bullet}]}\frac{h^ {|t|}}{\sigma(t)}D^2S(t),\\
 \tilde{p}(h) &=&p_2+p_3+\sum_{t\in [RT_{\circ}]}\frac{h^ {|t|}}{\sigma(t)}D^2S(t),
 \end{eqnarray}
respectively. 
\end{Lemma}

\begin{proof}


Uniqueness is trivial. Let us check that we have the right formal series. We only check 
equation (1). The other computation is similar.
\begin{eqnarray*}
\bar{x}(h) & = & x+h\nabla_{p_1}S(\bar{p},p_{3},x) \\
               & = & x+h\sum_{m\geq 0}\frac{1}{m!}\nabla_p^{(m+1)}S\bigg(\sum_{t\in [RT_{\circ}]}\frac{h^ {|t|}}{\sigma(t)}D^1S(t)%
               ,\dots,\sum_{t\in [RT_{\circ}]}\frac{h^ {|t|}}{\sigma(t)}D^1S(t)\bigg)\\
               & = & x+\sum_{m\geq 0}\sum_{t_1\in [RT_\circ]}\dots\sum_{t_m\in [RT_\circ]}\frac{h^ {1+|t_1|+\dots+|t_m|}}{m!\sigma(t_1)
               \dots\sigma(t_m)}\\
               &   & \times\nabla_p^{(m+1)}(D^1S(t_1),\dots,D^1S(t_m))\\
               & = & x +\sum_{m\geq 0}\sum_{t_1}\dots\sum_{t_m}\frac{h^{|t|}}{m!\sigma(t)}(\mu_1!\mu_2!\dots) D^1S(t),\quad \textrm{with }t=[t_1,\dots,t_m]_\bullet\\
               & = & x + \sum_{t\in[RT_\bullet]}\frac{h^{|t|}}{\sigma(t)}  D^1S(t)\\
\end{eqnarray*}

\end{proof}

We now insert these expansions into $M^1$ and $M^2$.
\begin{Lemma}
$$M^i(S)(h) = \sum_{t\in [RT]}
\frac{h^{|t|}}{\sigma(t)}
S^i(t)-\Big(\sum_{t\in [RT_\circ]} 
\frac{h^{|t|}}{\sigma(t)} D^iS(t) \Big) \Big(\sum_{t\in [RT_\bullet]}
\frac{h^{|t|}}{\sigma(t)} D^iS(t) \Big)$$ for $i= 1,2.$
\end{Lemma}

\begin{proof}
Let us do the proof for $M^1$. First we compute the different terms arising
in the formula for $M^1$ in terms of trees.
\begin{eqnarray*}
hS(p_1,p_2,\bar{x})  & = & h\sum_{m\geq 0}\frac{1}{m!}\nabla_x^{(m)}S\bigg(\sum_{t\in [RT_{\bullet}]}\frac{h^ {|t|}}{\sigma(t)}B(t)D^1S(t) ,\dots\\
                &   & ,\dots,\sum_{t\in [RT_{\bullet}]}\frac{h^ {|t|}}{\sigma(t)}B(t)D^1S(t)\bigg)\\
                & = & \sum_{m\geq 0}\sum_{t_1\in [RT_\bullet]}\dots\sum_{t_m\in [RT_\bullet]}\frac{h^{|t|}}{m!\sigma(t)}(\mu_1!\mu_2!\dots)\\
                &    & \times\nabla_x^{(m)}S(D^1S(t_1),\dots,D^1S(t_m)),\quad \textrm{with }t=[t_1,\dots,t_m]_\bullet\\ 
                & = & \sum_{t\in [RT_\bullet]}\frac{h^{|t|}}{\sigma(t)} S^1(t)
\end{eqnarray*}

By the same sort of computations we also get

\begin{eqnarray*}
h S(\bar{p},p_3,x) & = & \sum_{t\in [RT_\bullet]} \frac{h^{|t|}}{\sigma(t)} S^1(t)\\
h \nabla_xS(p_1,p_2,\bar{x}) & = & \sum_{t\in [RT_\circ]} \frac{h^{|t|}}{\sigma(t)} D^1S(t)\\
h \nabla_{p_1}S(\bar{p},p_3,x) & = &\sum_{t\in [RT_\bullet]} \frac{h^{|t|}}{\sigma(t)} D^1S(t)\\
\end{eqnarray*}

\end{proof}

The $M^i$'s are expressed as sums over topological rooted bipartite trees. We would like now to regroup
the terms of the formula in the previous Lemma. To do so we express all terms in terms of topological trees (no longer rooted).

\begin{Lemma}
Let $u\in[RT_\circ]$ and $v\in [RT_\bullet]$. Then $D^iS(u)D^iS(v) = S^i(u\circ v) = S^i(v\circ u)$.
\end{Lemma}

\begin{proof}
Prove it only for $i=1$. Suppose $u=[u_1,\dots,u_m]_\circ$, $v=[v_1,\dots,v_l]_\bullet$ then

\begin{eqnarray*}
D^1S(u)D^1S(v) & = & \nabla_x^{(m+1)}S(D^1S(u_1),\dots,D^1S(u_m)).D^1S(v)\\
               & = & \nabla_x^{(m+1)}S(D^1S(u_1),\dots,D^1S(u_m),D^1S(v))\\
               & = & S^1(u\circ v).\\
\end{eqnarray*}

\end{proof}

\begin{Lemma}\label{Lemma:sym}
Let $t=(V_t,E_t)\in T$. For all $v\in V_t$ let $t_v$ be the bipartite rooted tree $(V_t,E_t,v)\in RT$.
For $v\in V_t$ and $e=\{u,v\}\in E_t$ we have

\begin{eqnarray*}
 \frac{|sym(t)|}{|sym(t_v)|}  & = & |\{v'\in V_t /t_{v'}\textrm{is isomorphic to } t_v\}|\\
 \frac{|sym(t)|}{|sym(t_u)||sym(t_v)|}  & = & |\{e'\in E_t /t_{u'}\sqcup t_{v'}\textrm{is isomorphic to } t_u\sqcup t_v\}|
\end{eqnarray*}
\end{Lemma}

\begin{proof}
Consider the induced action of the symmetry group of the tree on the set of
vertices. Notice that two vertices $v$ and $w$ are in the same orbit iff
$t_v$ is isomorphic to $t_w$. Then the number of vertices of $t$ which lead to rooted
tree isomorphic to $t_v$ is exactly the cardinality of the orbit of $v$, which is exactly
$|sym(t)|$ divided by the cardinality of the isotropy subgroup which fixes $v$. But the latter
is $|sym(t_v)|$ by definition. We then get the first statement.
\par
For the second statement we have to consider the induced action on the edges and apply the
same type of argument. 
\end{proof}

\begin{Lemma}
Let $S\in C^\infty(B_2)[[h]]$.  The SGA equation for $S$ can be expressed
in terms of bipartite Cayley trees as
$$\sum_{t\in T}\frac{h^{|t|}}{|t|!}\big(S^1(t)-S^2(t)\big)=0.$$
\end{Lemma}

\begin{proof}
We have for $i=1,2$
\begin{eqnarray*}
M^i(S) & = & \sum_{t\in [RT]} \frac{h^{|t|}}{\sigma(t)}S^i(t)-%
             \sum_{u\in[RT_\circ]}\sum_{v\in[RT_\bullet]} \frac{h^{|u|+|v|}}{\sigma(u)\sigma(v)}D^iS(u)D^iS(v)\\
        & = & \sum_{\bar{t}\in [T]} h^{|\bar{t}|}S^i(\bar{t})%
                \Big\{ \sum_{t\in\bar{t}} \frac{1}{|sym(t)|}-\sum_{\substack{u\in [RT_\bullet],v\in [RT_\circ]\\ u\circ v \in\bar{t}}}%
                   \frac{1}{|sym(u)||sym(v)|}   \Big\}\\
        & = & \sum_{t\in T}\frac{h^{|t|}}{|t|!}S^i(t)\Big\{%
        \sum_{v\in V_t} \frac{|sym(t)|}{|sym(t_v)|}\frac{1}{k(t,v)}\\
        &    & -\sum_{e = \{u,v\}\in E_t} \frac{|sym(t)|}{|sym(t_u)||sym(t_v)|}\frac{1}{l(t,e)} \Big\}
\end{eqnarray*}
where $k(t,v) = |\{v'\in V_t /t_{v'}\textrm{is isomorphic to } t_v\}|$ 
and $l(t,e) = |\{e'\in E_t /t_{u'}\sqcup t_{v'}\textrm{is isomorphic to } t_u\sqcup t_v\}|$.
Using Lemma \ref{Lemma:sym} and the fact that for a tree the difference between the number of vertices
and the number of edges is equal to 1 we get the desired result.
\end{proof}

Using now the fact that $S$ is a formal series we immediately get Proposition \ref{pertSGA}.

\section{Geometry of Kontsevich trees}\label{Kont}

\label{sec:Konttrees}
In this section we present a diagrammatical notation introduced by Kontsevich which 
allows us to write an explicit solution of the SGA equation. 

\subsection{Basic Definitions}
      \begin{Definition}
      \begin{enumerate}

      \item A \textbf{Kontsevich  graph $\Gamma$ of type $(n,m)$} is a directed graph $\Gamma = 
      (E_{\Gamma},V_{\Gamma})$ which has the following properties:
        \begin{itemize}
                \item it possesses two types of vertices $V_{\Gamma}=V_{\Gamma}^a\ \sqcup V_{\Gamma}^g$,
                the aerial vertices $V_{\Gamma}^a = \{1,\dots,n\}$ and the ground vertices $V_{\Gamma}^g = \{\bar{1},\dots,\bar{m}\}$ .
                \item each aerial vertex possesses exactly two ordered edges starting from it. 
                The edge set can be described as $E_{\Gamma} = \{(k,\gamma^i(k)),\quad k = 1,
                \dots,n,\quad i = 1,2\}$ where $\gamma^i:V_{\Gamma}^a\rightarrow V_{\Gamma}$. Sometimes
                one denotes the two edges of a vertex $k$ by $e^1_k$ and $e^2_k$.
                \item For each aerial vertex $v$ we do not allow small loops (i.e., that $\gamma^i(v) = v$)  and
                 double edges (i.e., that $\gamma^1(v) = \gamma^2(v)$).
        \end{itemize}
        We denote the set of Kontsevich graphs of type $(n,m)$ by $G_{n,m}$.
        If $\Gamma\in G_{n,m}$ then we set $|\Gamma| := n$.

        \item Let $A \in V_\Gamma$. We call $\Gamma_{/A}$ the \textbf{restriction} of $\Gamma$ to $A$. It is the 
        graph with vertex set $A$ and edges $E_\Gamma \cap A \times A$. We call $\Gamma_{(A)}$ the \textbf{contraction}
        of $\Gamma$ to $A$. It is the graph with vertex set $(V_\Gamma\backslash A)\sqcup  \{*\}$ (the vertices  of $A$ are contracted
        to a single vertex $*$) and edges $(i,j)\in E_\Gamma$ where $i$ is replaced by the new vertex $*$ in $\Gamma_{(A)}$ if $i\in A$
        and the same for $j$ (simple loops are deleted). Note that the resulting graphs might not be Kontsevich graphs.

      \item We denote  by $\Delta(\Gamma) = (V_{\Gamma}^a,E_{\Gamma}^a)$ the restriction
      of $\Gamma\in G_{n,m}$ to the aerial vertices. Sometimes we write  $E_{\Gamma}^g = E_\Gamma\backslash E_\Gamma^a$. We 
      say that a Kontsevich graph $\Gamma$ is \textbf{connected} if $\Delta(\Gamma)$ is connected in the usual sense.
      We say that a connected Kontsevich graph $\Gamma$ is a tree if $\Delta(\Gamma)$ is a tree(i.e., a graph without cycle).
        Denote by $C_{n,m}$ the set of connected Kontsevich graph of type $(n,m)$ and 
      by $T_{n,m}$ the set of Kontsevich trees of type $(n,m)$.

      \end{enumerate}

      \end{Definition}

Given a Poisson structure $\alpha$ on $\mathbb R^d$ one can associate to each graph $\Gamma\in G_{n,m}$ an
$m$-multidifferential operator on $C^\infty(\mathbb R^d)$. The general formula is the following
$$
B_{\Gamma}(f_1\dots,f_m) := \sum_{I:E_{\Gamma}\rightarrow\{1,\dots,d\}} \big[\prod_{k\in V_\Gamma^a}(\prod_{\substack{e\in E_{\Gamma} \\ e=(*,k)}}\partial_{I(e)})\alpha^{I(e_k^1) I(e_k^2)}\big]\times  
                  \prod_{i\in V_\Gamma^g}\big(\prod_{\substack{e\in E_{\Gamma} \\ e=(*,i)}}\partial_{I(e)}\big)f_i
$$
We call $\hat B_\Gamma$  the \textbf{symbol} of $B_\Gamma$. It can be defined by the
formula
$$
B_\Gamma (e^{p_1x},\dots,e^{p_mx}) = \hat B_\Gamma(p_1,\dots,p_m,x)e^{(p_1+\dots+p_m)x}$$

\vspace{0.3cm}
\begin{Example}
Take the graph
$$\Gamma = \kont{K3-2}$$ then we have $$ \hat B_\Gamma(p_1,p_2,x) = \sum_{1\leq i,j,k,l,m,n \leq d}\alpha^{ij}(x) \partial_n \partial_j \alpha^{kl}(x) \alpha^{mn}(x) p^1_k  p^2_i p^2_l  p^2_m .$$
\end{Example}

Associated to each Kontsevich graph $\Gamma\in G_{n,m}$ there is also a number, the \textbf{Kontsevich weight} $W_\Gamma$.
In these notes we only need to define these weights for graphs of type $(n,2)$. The generalization is however straightforward.
We do this in several steps.

        \begin{enumerate}
                \item{
                Take a Kontsevich graph $\Gamma\in G_{n,2}$ and identify its
                vertices $1,\dots,n\in V_{\Gamma}$ with $n$ complex numbers
                $z_1,\dots,z_n$ lying in the upper half complex plane
                $\mathcal{H}=\{z\in \mathbb{C}\quad / \quad Im(z) > 0\}$(we require that $z_i\neq z_j$ if 
                $i\neq j$). Identify further $\bar 1$ and $\bar 2$ with $0$ and
                $1$ in $\mathbb R$.
                }

                \item{
                Consider now the hyperbolic metric on $\mathcal{H}$.
                The geodesic joining two points $p,q\in\mathcal{H}$ is
                in this metric either the half circle intersecting orthogonally
                the real line and passing through $p$ and $q$ or the 
                 line orthogonal to the real line passing through $p$ and $q$.
                We can now associate the oriented
                edges $e_k^i=(k,\gamma^i(k))$ to the oriented geodesics
                joining $z_k$ and $z_{\gamma^i(k)}$. We call such an
                embedding of $\Gamma$ a $\it{configuration}$ of $\Gamma$.
                We can then identify the configuration space of a Kontsevich 
                graph $\Gamma$ with $\mathcal{H}^n\backslash D^n$
                where $\mathcal{H}^n$ is $n$ times the Cartesian 
                product of $\mathcal{H}$ and
                $$D^n:=\{(z_1,\dots,z_n)\in\mathcal{H}^n\quad/\quad
                \exists i,j\quad i\neq j\quad \textrm{and}\quad z_i=z_j\}.$$
                Notice that $\mathcal{H}^n\backslash D^n$ is a real
                non-compact manifold of dimension $2n$. We can however compactify it into
                a compact manifold with corners $\overline{\mathcal H^n\backslash D^n}$ such that the 
                open stratum is exactly $\mathcal H^n\backslash D^n$.
                }

                \item{
                For each edge $e_k^i=(k,\gamma^i(k))$ we can define
                an ``angle function" on $\mathcal H^n\backslash D^n$ by
                $\psi^i_{z_k}(z_1,\dots,z_n):=\phi^h(z_k,z_{\gamma^i(k)})$
                where $\phi^h(z_k,z_{\gamma^i(k)})$ is the oriented hyperbolic
                angle between the geodesic joining $z_k$ and $\infty$ and
                the geodesic joining $z_k$ and $z_{\gamma^i(k)}$. So 
                $\phi^h(p,q)=\arg \big( \frac{q-p}{q-\bar{p}}\big).$
                }

                \item{
                We can now consider the $1$-forms 
                $d\psi^i_{z_k}\in \Omega^1(\mathcal{H}^n\backslash D^n)$
                which can be extended on the compactified space. Then the Kontsevich
                weight of $\Gamma$ is defined by

                $$W_{\Gamma}:=\frac{1}{(2\pi)^{2n}}\int_{\overline{\mathcal{H}^n\backslash D^n}}
                \bigwedge_{i=1}^n(d\psi^1_{z_k}\wedge d\psi^2_{z_k}).$$



                }
        \end{enumerate}

        Further explanations about these operators and weights can be found in \cite{kontsevich1997}.
        However we still need a Lemma which is also proven in (or follows directly from) 
        \cite{kontsevich1997}.

\begin{Definition}
Let $\Gamma\in G_{n,3}$. We denote by $\sub(\Gamma)_{\{\bar1,\bar2\}}$ the set of the subset $S$ of $V^a_\Gamma$
such that $\Gamma_{/\{\bar1,\bar2\}\sqcup S}$ and $\Gamma_{(\{\bar1,\bar2\}\sqcup S)}$ are still
Kontsevich graphs of type $(n,2)$. We define similarly $\sub(\Gamma)_{\{\bar 2,\bar 3\}}$
\end{Definition}

\begin{Lemma}\label{vanishing}
\[
\sum_{\Gamma\in G_{n,3}}\Big(
\sum_{S\in \sub(\Gamma)_{\{\bar1,\bar2\}} }W_{\Gamma_{/\{\bar1,\bar2\}\sqcup S}}W_{\Gamma_{(\{\bar1,\bar2\}\sqcup S)}}
-
\sum_{S\in \sub(\Gamma)_{\{\bar2,\bar3\}} }W_{\Gamma_{/\{\bar2,\bar3\}\sqcup S}}W_{\Gamma_{(\{\bar2,\bar3\}\sqcup S)}}
\Big)\hat B_\Gamma = 0
\]

\end{Lemma}

\subsection{Factorization into connected components  of graphs of type $(n,2)$ }

We describe here a procedure which allows us to decompose a graph of type $(n,2)$
into $l$ graphs $\Gamma_1,\dots,\Gamma_l$ of the same type, its connected components
in a slight unusual sense. Take $\Gamma\in G_{n,2}$. Then

      \begin{enumerate}
         \item Consider the usual connected components  of 
                 $\Delta(\Gamma)$. We can number them in a unique way
                 using the following rule:  Let $\Delta^i (\Gamma)$, $\Delta^j (\Gamma)$ be two
                 connected components of $\Delta(\Gamma)$. We impose that $i<j $ iff $ \min\{V_{\Delta^i (\Gamma)}\} <
                 \min\{V_{\Delta^j (\Gamma)}\}$

         \item For each connected component $\Delta^i (\Gamma)$ of $\Delta(\Gamma)$ we can 
                reconstruct a Kontsevich graph which we denote by $\Gamma_i$:
                \begin{enumerate}
                    \item To begin with, add to each $\Delta^i (\Gamma)$ 
                    the vertices and edges that we removed considering $\Delta(\Gamma)$.
                           Let $\hat{\Gamma}_i$ be this graph.
                          
                   \item Relabel the vertices of $\hat \Gamma_i$ by $1,2,\dots,|\Delta^i(\Gamma)|$ preserving
                   the relative order of the vertices of $\Delta^i(\Gamma)$. One gets a new Kontsevich graph
                   $\Gamma_i$.
                        
                 \end{enumerate}



      \end{enumerate}

\begin{Definition}
          \begin{enumerate}
                \item Let $\Gamma \in G_{n,2}$. We call the $\Gamma_i$'s as constructed above the
                \textbf{connected factors} of $\Gamma$.  Because of the numbering of the $\Delta^i(\Gamma)$ the connected 
                factors of a Kontsevich graph $\Gamma$ are uniquely numbered.
                The connected factors of $\Gamma$ are connected Kontsevich graphs.

     \item We denote by $G_{n,2}(n_1,\dots,n_k)$ the graphs $\Gamma$ of $G_{n,2}$ which
        have $k$ connected factors and such that the $i^{th}$ connected factors $\Gamma_i$
         is a Kontsevich  graph of order $n_i$.

         \item We call the \textbf{factorization map} the map $D$ defined by
        $ D(\Gamma) = (\Gamma_1, \dots, \Gamma_k)$ where the $\Gamma_i$ are the connected factor of $\Gamma$.
      \end{enumerate}

 \end{Definition}
      Similar considerations about connected Kontsevich graphs and connected factorization
      can be found in \cite{kathotia1998}. In particular one can find the following Lemma:

\begin{Lemma}[Factorization Lemma]\label{lem:fact}
Let $\Gamma \in G_{n,2}$ and $D(\Gamma)=(\Gamma_{1},\dots,\Gamma_{k})$
its connected factorization. Then we have
\begin{enumerate}
\item $\mathcal{W}_{\Gamma}=\mathcal{W}_{\Gamma_{1}} \dots \mathcal{W}_{\Gamma_{k}}$
\item $\hat B_\Gamma = \hat{B}_{\Gamma_{1}} \dots \hat{B}_{\Gamma_{k}}$.
\end{enumerate}
\end{Lemma}

\subsection{Number of graphs leading to the same connected factorization}

We are looking for the number of graphs of $G_{n,2}$ which lead to
the same connected factorization. This number plays a crucial  role while
proving the existence of the generating function.

        It is clear that $D(\Gamma) = D(\Gamma')$ only if 
        $\Gamma, \Gamma'  \in G_{n,2}(n_1,\dots,n_k)$ for some
        $n_1,\dots,n_k$.
        Therefore the problem of counting the number of Kontsevich graphs of type
        $(n,2)$ that lead to the same factorization can be stated in the following terms:

        \vspace{10pt}
        \emph{Given} $(\Gamma_1, \dots, \Gamma_k) \in C_{n_1,2} \times \dots \times C_{n_k,2},$ 
        \emph{what is the number of elements of} $D^{-1}(\Gamma_1, \dots, \Gamma_k)$? 
        \vspace{10pt}

        The answer is contained in the following remarks.
\par

     Notice that the permutation group $S_n$ acts on $G_{n,2}$ by permuting
     the aerial vertices.

\par
        Let $\Gamma \in G_{n,2}(n_1,\dots,n_k)$. All the graphs 
        $\Gamma'  \in G_{n,2}(n_1,\dots,n_k)$ which give the same
        connected factorization as $\Gamma$ are generated by a subset of $S_n$, i.e.,
        $$\forall \Gamma'  \in G_{n,2}(n_1,\dots,n_k) \textrm{ s.t. } 
        D(\Gamma) = D(\Gamma') \quad \exists \sigma \in P \textrm{ s.t. }
        \sigma \Gamma = \Gamma'.$$

        This subset $P \subset S_n$ is defined by the constraints:
        \begin{enumerate}
                \item The permutation must preserve the relative order
                of the vertices of  $V_{\Gamma_i}$.

                \item Consider the set of the minimum vertex of each $V_{\Gamma_i}$. 
                The permutation must preserve the relative order of this set.
        \end{enumerate}

        It remains then to count the number of such permutations.
        The second constraint restricts the number of allowed permutations 
        to $\frac{n!}{k!}$. The first further restricts to
        $\frac{n!}{k!n_1!\dots n_k!}$. Thus $$|D^{-1}(\Gamma_1, \dots, \Gamma_k)| = \frac{n!}{k!n_1!\dots n_k!}.$$
        As this number reappears in another context let us denote it by $d(n_1,\dots,n_k)$
and call it the \textbf{decomposition coefficient}.


\subsection{Contraction-Restriction decomposition of trees of type $(n,3)$}

 Here begin some new considerations about Kontsevich graphs. We will see that in each
Kontsevich tree of type $(n,3)$ lies, hidden, two Cayley trees which encode the 
contraction and restriction of the tree leading to Kontsevich trees. These two Cayley trees allow us to make a link between the
perturbative SGA equation which is expressed in terms of Cayley trees and the proposed solution
expressed in terms of Kontsevich trees. The main results of this section are then
summarized in definition \ref{def:CRpattern} and Proposition \ref{konttrees}. 
But let us begin first to establish a few little facts necessary to make any statement.

      \begin{Lemma}
      Let $\Gamma\in T_{n,m}$ then
      \begin{enumerate}
      \item $|E_{\Gamma}^a|=n-1$
      \item $|E_{\Gamma}^g|=n+1$
      \end{enumerate}
      \end{Lemma}

\begin{proof}
 For the  first assertion one notices that $\Delta(\Gamma)$,
which has $n$ vertices, is  connected, so there are at least $n-1$ edges connecting these vertices. Now, if we
add an edge, we create a cycle which contradicts the fact that $\Delta(\Gamma)$ is  a tree.
The second assertion follows from the identity $|E^a_\Gamma|+|E^g_\Gamma|=2n$.
\end{proof}

\begin{Corollary}\label{cor:tree1}
There is no Kontsevich tree of type $(n,1)$ (i.e. $T_{n,1} = \emptyset$).

\end{Corollary}

\begin{proof}

        As $|E_\Gamma^g|=n+1$ and $|V_\Gamma^a|=n$, one aerial vertex has its two
        edges landing at the only ground vertex and we do not allow double edges.
\end{proof}

\begin{Corollary}
Suppose $\Gamma\in T_{n,2}$. Then  $E^g_\Gamma$ has at
least one edge landing at $\bar1$ and one edge landing at $\bar{2}$.
\end{Corollary}

\begin{proof}
Without loss of generality, suppose that all edges of $E^g_\Gamma$ land at  $\bar1$ then $\Gamma_{/V^a_\Gamma\sqcup \bar1}\in T_{n,1}
=\emptyset$.
\end{proof}

\begin{Corollary}\label{cor:tree3}
Suppose $\Gamma\in T_{n,2}.$ There is at least one $v\in V_{\Gamma}^a$ such that $\gamma^1(v) =\bar{1}$, 
$\gamma^2(v) =\bar{2}$.
\end{Corollary}

\begin{proof}
As $|E_\Gamma^g|=n+1$ and $|V_\Gamma^a|=n$, there is one aerial vertex whose both edges are
ground edges. Those two edges can not land at  the same ground vertex as we prevent double edges.
\end{proof}

\begin{Definition}
\begin{enumerate}
\item Let be $\Gamma\in G_{m,n}$. One defines the following transitive relation among the vertices
of $\Gamma$: $v<w$ iff there exists $a_1,\dots,a_k\in V_\Gamma$ such that 
$$(w,a_1),\dots,(a_i,a_{i+i}),\dots,(a_k,v)\in E_\Gamma.$$

\item Let be $\Gamma\in G_{n,m}$. Let us denote by  
$$\Star_{in}(v) := \{w\in V_\Gamma\quad s.t. \quad v<w\}$$
$$\Star_{out}(v) := \{w\in V_\Gamma\quad s.t. \quad w<v\}.$$
\end{enumerate}
\end{Definition}

\begin{Lemma}\label{lem:boulesblanches}
Let $\Gamma\in T_{n,3}.$ Denote $N_{\bar{1}} := \Star_{in}(\bar{1})$, $B_{\bar 1}:=V_\Gamma^a\backslash N_{\bar{1}}$
and $\Gamma^1_{B_{\bar{1}}},\dots,\Gamma^l_{B_{\bar{1}}}$ the connected factors of $\Gamma_{/\{\bar{2},\bar{3}\}\sqcup B_{\bar{1}}}$.
Then the $\Gamma_{B_{\bar{1}}}^i$'s are Kontsevich trees with two ground vertices (provided that $B_{\bar{1}} \neq \emptyset$).
The same statement holds if we replace $B_{\bar{1}}$ by $B_{\bar{3}}$ and make the restriction around $\{\bar{1},\bar{2}\}\sqcup B_{\bar{3}}$.
\end{Lemma}

\begin{proof}
Take $\Gamma^i_{B_{\bar1}}$. As there are no edges $(v,w)$ starting from $B_{\bar1}$ and landing at $N_{\bar1}\sqcup \bar1$ (otherwise 
$v>w>\bar1 \Rightarrow v\in N_{\bar1}$), all the vertices of $B_{\bar1}$ conserve their two edges when passing
to the restriction $\Gamma_{/\{\bar2,\bar3\}\sqcup B_{\bar1}}$.
It remains to be shown that all the edges $E_{\Gamma^i_{B_{\bar1}}}^g $ are not landing exclusively at one
of $\bar1$ or $\bar2$. But corollary \ref{cor:tree3} prevents  this phenomenon from happening.
\end{proof}

\textbf{\emph{Trivial little facts:}}
\par
We define for convenience $B^i_{\bar{1}}:=V_{\Gamma^i_{B_{\bar 1}}}^a$, $i=1,\dots,l$ and 
$N^j_{\bar{1}}:=V_{\Gamma^j_{N_{\bar 1}}}^a$, $j=1,\dots,k$ where $\Gamma^j_{N_1}$ are the connected
factors of $\Gamma_{(\{\bar{2},\bar{3}\}\sqcup B_{\bar{1}})}$. We see that:
\begin{enumerate}
\item There is at most one edge from $\Gamma$ starting from one $N^j_{\bar{1}}$ to a $B^i_{\bar{1}}$( otherwise one introduces a cycle).
\item There is no edge from an $N^j_{\bar{1}}$ to another $N^i_{\bar{1}}$ (they are connected factors).
\item There is no edge from a $B^j_{\bar{1}}$ to another $B^i_{\bar{1}}$ (they are connected factors).
\item There is no edge from  a $B^j_{\bar{1}}$ to a $N^i_{\bar{1}}$ (otherwise one vertex of  $B^j_{\bar{1}}$ should
be in $N^i_{\bar{1}}$).
\end{enumerate}

\begin{Corollary}[Contraction/Restriction trees]
Let $\Gamma\in T_{n,3}$. We can make the following construction:

\begin{itemize}
        \item identifying each $N_{\bar{1}}^j$, $j=1,\dots,k$ and $B_{\bar{1}}^i$, $i=1,\dots,l$ with
        respectively black vertex and white vertex,
        \item putting an edge between black vertex and white
        vertex iff there is one edge between the corresponding sets $N_{\bar{1}}^j$ and $B_{\bar{1}}^i$,
        \item labelling the black and white vertices such that $i<j$ iff the minimum of the set corresponding to $i$ is
        inferior to the minimum of the set corresponding to $j$,
\end{itemize}
we get a Cayley tree $t_\Gamma^2\in T$.
This tree $t^2_\Gamma$ is  called the second contraction/restriction tree of $\Gamma$.
If we start the construction from $B_{\bar{3}}$ and $N_{\bar{3}}$ we get $t^1_\Gamma$, the 
first contraction/restriction tree of $\Gamma$.

\end{Corollary}

\begin{Example}
The following graph $\Gamma$ illustrates these phenomenon.

\vspace{5cm}
\kont{WB-1}
\vspace{0.5cm}

For this graph we have that the two contraction/restriction trees are

$$t_\Gamma^1 = \bullet \quad\textrm{and}\quad t_\Gamma^2 = \rtree{RT-4}.$$
\end{Example}

\begin{Lemma}\label{lem:boulesnoires}
Let $\Gamma\in T_{n,3}$. Denote $N_{\bar{1}} := \Star_{in}(\bar{1})$, $B_1:=V_\Gamma^a\backslash N_{\bar{1}}$ and
 $\Gamma^1_{N_{\bar{1}}},\dots,\Gamma^1_{N_{\bar{k}}}$ the connected factor of $\Gamma_{(\{\bar{2},\bar{3}\}\sqcup B_{\bar{1}})}$.

Then the $\Gamma_{N_{\bar{1}}}^i$'s are Kontsevich trees with two ground vertices (provided that $B_{\bar{1}} \neq \emptyset$).
The same statement holds if we replace $B_{\bar{1}}$ by $B_{\bar{3}}$ and make the contraction around $\{\bar{1},\bar{2}\}\sqcup B_{\bar{3}}$
\end{Lemma}

\begin{proof}
From the vertices in $N^i_{\bar1}:=V_{\Gamma_{N_{\bar1}}^i}$, there is 
at least one edge landing at $\bar 1$ and at most one landing at each $B_{\bar1}^\mu$.
The only bad thing that can happen is then that
there is $v\in N_{\bar1}^i$ such that $\gamma^1(v)\in B_{\bar1}^\mu$ and $\gamma^1(v)\in B_{\bar1}^\nu$.
But then $v$ has no any edge left starting from it. Which implies that $\bar1\notin \Star_{out}(v)$.
\end{proof}

\begin{Definition}\label{def:CRpattern}
Let $\Gamma\in T_{n,3}$. We define the contraction/restriction decomposition maps
$$P^i(\Gamma)=(t_\Gamma^i,\Gamma_1,\dots,\Gamma_m),\quad i=1,2$$
where $t_\Gamma^i\in T$ is the $i^{th}$ contraction/restriction-tree of $\Gamma$ and
the $\Gamma^j$ are the connected factor of the contraction and the restriction of $\Gamma$ 
around $\{\bar1,\bar2\}\sqcup B_{\bar3}$ for $i=1$ and around $\{\bar2,\bar3\}\sqcup B_{\bar1}$ for $i=2$.
We index these connected factors with the usual convention, that is $k<l$ if the minimum of the aerial vertices of 
$\Gamma_k$ is less than the minimum of the aerial vertices of $\Gamma_l$.

\par
We denote by $T^i_{n,3}(t,\Gamma_1,\dots,\Gamma_{|t|})$ the subset of $T_{n,3}$ such that $P^i(\Gamma) = (t,\Gamma_1,
\dots,\Gamma_{|t|})$ for $i=1,2$.

\end{Definition}

\begin{Example}
For the previous graph $\Gamma$ we get

\begin{eqnarray*}
P^1(\Gamma) & = & \Bigg(\bullet_1,\kont{K8-1}\quad\Bigg)\\
&&\\
&&\\
P^2(\Gamma)  & = &  \Bigg(\rtree{RT-4}\quad, \kont{K2-1},\kont{K1-1},\kont{K3-7},\kont{K2-2}\Bigg)
\end{eqnarray*}
\end{Example}

\begin{Proposition}\label{konttrees}
Let $\Gamma\in T_{n,3}$. Then in the notation used above we have

\begin{enumerate}
        \item Let $\Gamma\in T^1_{n,3}(t;\Gamma_1,\dots,\Gamma_{|t|})$ then
        $$W_{\Gamma_1}\dots W_{\Gamma_{|t|}}=W_{\Gamma /B_{\bar 3}\cap \{\bar 1,\bar 2\}}W_{\Gamma (B_{\bar 3}\cap \{\bar 1,\bar 2\})}$$
         \par
         Let $\Gamma\in T^2_{n,3}(t;\Gamma_1,\dots,\Gamma_{|t|})$ then
        $$W_{\Gamma_1}\dots W_{\Gamma_{|t|}}=W_{\Gamma /B_{\bar 1}\cap \{\bar 2,\bar 3\}}W_{\Gamma (B_{\bar 1}\cap \{\bar 2,\bar 3\})}$$

        \item We have the following equations for the Kontsevich weights 
        $$\sum_{\Gamma\in T_{n,3}}\Big(W_{\Gamma_{(\{\bar1,\bar2\}\sqcup B_{\bar3})}}W_{\Gamma_{/\{\bar1,\bar2\}\sqcup B_{\bar3}}}
        -W_{\Gamma_{(\{\bar2,\bar3\}\sqcup B_{\bar1})}}W_{\Gamma_{/\{\bar2,\bar3\}\sqcup B_{\bar1}}}\Big)\hat B_\Gamma =0$$

        \item The following relates Cayley trees and Kontsevich trees, for all $t\in T$ we have
        $$C^i_t(\hat B_{\Gamma_1},\dots,\hat B_{\Gamma_{|t|}}) = d(n_1,\dots,n_{|t|})\sum_{\Gamma\in T^i_{n,3}(t;\Gamma_1,\dots,\Gamma_{|t|})}\hat B_\Gamma$$
\end{enumerate}

\end{Proposition}

\begin{proof}
\begin{enumerate}

\item[(1)] is trivial.\par
\item[(2)] is a consequence of Lemma \ref{vanishing} once one has proved that $\sub(\Gamma)_{\{\bar 1,\bar 2\}}= \{B_{\bar 3}\}$
and $\sub(\Gamma)_{\{\bar 2,\bar 3\}}= \{B_{\bar 1}\}$. By the Lemmas \ref{lem:boulesblanches} and \ref{lem:boulesnoires} one has already
that $B_{\bar 1}\in\sub(\Gamma)_{\{\bar 2,\bar 3\}}$ and $B_{\bar 3}\in\sub(\Gamma)_{\{\bar 2,\bar 3\}}$. It remains to check that they are the only ones.  Let us prove that only for $B_{\bar1}$.
\par
Suppose  there is another subset $K\subset V_\Gamma^a$ such that 
$\Gamma_{(\{\bar2,\bar3\}\sqcup K)}$ and $\Gamma_{/\{\bar2,\bar3\}\sqcup K}$ are Kontsevich
trees. This implies that in the process of
\begin{enumerate}
        \item restriction around $\{\bar2,\bar3\}\sqcup K$, one should not loose an edge
        \item contraction around $\{\bar2,\bar3\}\sqcup K$, one should not end up with a double edge 
\end{enumerate}

(A) Suppose that $K\cap N_{\bar1} \neq \emptyset$. Take  $v\in K\cap N_{\bar1}$ then 
$\Star_{out}(v)$ is a subset of $K$ otherwise we loose an edge when doing the restriction around 
$\{\bar2,\bar3\}\sqcup K$. But $\bar1 \in \Star_{out}(v)$ which implies that $\bar1\in K$
otherwise we loose an edge when doing the restriction. Contradiction with $K\subset V^a_\Gamma$.
\par

(B)By (A) we have that $K\subset B_{\bar1}$. Suppose that $K$ is strictly contained in $B_{\bar1}$.
Then $(\Gamma_{/\{\bar2,\bar3\}\sqcup B_{\bar1}})_{(K\sqcup \{\bar2,\bar3\})}$ is a subgraph of $\Gamma_{(K\sqcup \{\bar2,\bar3\})}$.
But as there are no edge starting from $B_{\bar1}$ and landing at $\bar1$,
$(\Gamma_{/\{\bar2,\bar3\}\sqcup B_{\bar1}})_{(K\sqcup \{\bar2,\bar3\})}$ is a Kontsevich tree with only one ground
vertex which implies that it is not a Kontsevich tree. Contradiction.

\item[(3)]First remark that $\sum_{\Gamma\in T^i_{n,3}(t,\Gamma_1,\dots,\Gamma_{|t|})} 
B_{\Gamma}= d(n_1,\dots,n_{|t|})\sum_{\Gamma\in A} B_{\Gamma} $
where $A$ is the subset of trees $\Gamma\in T^i_{n,3}(t,\Gamma_1,\dots,\Gamma_{|t|})$ such that 
all vertices in $V_\Gamma$ corresponding to these of $V_{\Gamma_i}$ are less than these corresponding
to $V_{\Gamma_j}$ if $i<j$.
It is clear that letting act all the permutations of $S_n$ which preserve the relative order of the
minimal vertex of each $V_{\Gamma_i}$ and the relative order of the vertices in $V_{\Gamma_i}$
we get all trees of $T^i_{n,3}(t,\Gamma_1,\dots,\Gamma_{|t|})$. We have already
counted the number of such permutations it is exactly the decomposition coefficient $d(n_1,\dots,n_{|t|})$.

 The identity $\sum_{\Gamma\in A} B_{\Gamma} =  C^i_t(\hat B_{\Gamma_1},\dots,\hat B_{\Gamma_k})$ follows
 from the Leibniz rule.

\end{enumerate}
\end{proof}

\section{Proof of Theorem \ref{th:main}}

\label{proof}
Let us restate the main Theorem. 

\begin{Th}
Given a Poisson structure $\alpha$ on $\mathbb R^d$
there exists a unique natural deformation of the trivial generating function such that
the first order is precisely $\alpha$. Moreover we have an explicit formula for this deformation
$$S_h(p_1,p_2,x) = x(p_1+p_2)+\sum_{n=1}^\infty \frac{h^n}{n!}\sum_{\Gamma\in T_{n,2}}W_\Gamma \hat B_\Gamma(p_1,p_2,x)$$
where $T_{n,2}$ is the set of Kontsevich trees of type $(n,2)$, $W_\Gamma$ is the
Kontsevich weight of $\Gamma$ and $\hat B_\Gamma$ is the symbol of the bidifferential
operator $B_\Gamma$ associated to $\Gamma$.
\end{Th}

\begin{proof}
\textbf{Existence of the solution.}\par

Let us verify that the proposed solution satisfies the perturbative version of the SGA equation.
Denote 
$$M^i_n(S) = \sum_{\substack{t\in T\\ |t|\leq n}}\frac{1}{|t|!}\sum_{\substack{n_1+\dots+n_{|t|}=n\\n_i\geq 1}}C_t^i(S_{n_1},\dots,S_{n_{|t|}})
$$
Let us compute $M^1_n(S)$ for the proposed solution
\begin{eqnarray*}
M^1(S)_n & = & \sum_{\substack{t\in T\\ |t|\leq n}} \frac{1}{|t|!}\sum_{\substack{n_1+\dots+n_{|t|}=n\\n_i\geq 1}} \sum_{\substack{\Gamma_i\in T_{n_i,3}\\ i=1,\dots,|t|}}\frac{W_{\Gamma_1}\dots W_{\Gamma_{|t|}}}{n_1!\dots n_{|t|}!} C_t^1\Big(\hat B_{\Gamma_1},\dots,\hat B_{\Gamma_{|t|}}\Big)\\
         & = & \sum_{\substack{t\in T\\ |t|\leq n}} \sum_{\substack{n_1+\dots+n_{|t|}=n\\n_i\geq 1}} \sum_{\substack{\Gamma_i\in T_{n_i,3}\\ i=1,\dots,|t|}}\frac{W_{\Gamma_1}\dots W_{\Gamma_{|t|}}}{(n_1+\dots +n_{|t|})!}\sum_{\Gamma\in T^1_{n,3}(t;\Gamma_1,\dots,\Gamma_{|t|})}B_\Gamma \\
        & = & \sum_{\substack{t\in T\\ |t|\leq n}}\frac{1}{n!} \sum_{\substack{n_1+\dots+n_{|t|}=n\\n_i\geq 1}} \sum_{\substack{\Gamma_i\in T_{n_i,3}\\ i=1,\dots,|t|}}\sum_{\Gamma\in T^1_{n,3}(t,\Gamma_1,\dots,\Gamma_{|t|})} W_{\Gamma_{(\{\bar1,\bar2\}\sqcup B_{\bar3})}}W_{\Gamma_{/\{\bar1,\bar2\}\sqcup B_{\bar3}}} B_{\Gamma}\\
        & = & \sum_{n\geq 1}\frac{1}{n!}\sum_{\Gamma\in T_{n,3}}%
        W_{\Gamma_{(\{\bar1,\bar2\}\sqcup B_{\bar3})}}W_{\Gamma_{/\{\bar1,\bar2\}\sqcup B_{\bar3}}} B_{\Gamma}
\end{eqnarray*}
which implies by Proposition \ref{konttrees} that $M^1_n(S)-M^2_n(S)=0$ for all $n>0$.

\vspace{1cm}
\textbf{Uniqueness of the solution.}\par
We have seen that the perturbative SGA equations could be put at each order into the
form $dS_m+H_m(S_{m-1},\dots,S_1)=0$ where the differential $d$ may be identified with
the Hochschild differential on symbols. 

Let $S$ and $S'$ be two generating functions. By definition we have that $S_{1}=S'_{1}=\alpha$.
Now suppose that $S$ and $S'$ are equal up to order $m-1$ (i.e., $S_{k}=S'_{k}$, $k\leq m-1$). Thus
$K_{m}:=S_{m}-S'_{m}\in C^\infty(B_2)$ satisfies the following equation
$$
dK_{m}=H_{m}(S_{1},\dots,S_{m-1})-H_{m}(S'_{1},\dots,S'_{m-1})=0.
$$

As $H^ {2}(C^\infty(B_\bullet),d) = V^2(\mathbb R^d)$(bivector fields over $\mathbb R^d$)
we have that $K_m$ can be written as $ K_m = dk_m + \omega$
where $k_m$ is a $1$-cochain and $\omega$ is a bivector field.
Because of the homogeneity of $K_m$ in the $p$'s we have that 
$\omega$ vanishes. 
\par
\textbf{Claim:} $k_m(p):= \frac{-1}{m+1}K^1_m(p,p)$ is a primitive of $K_m$, i.e., $dk_m = K_m$.
\par
This claim prove the uniqueness because by assumption we have $K^1_m(p,p)=0$ which means that $k_m=0$
and thus $dk_m = K_m=0$. As for the claim, suppose that $K_m(p_1,p_2) = \sum_{|I|+|J|=m+1}K_m^{I,J}
\frac{p_1^Ip_2^J}{I!J!}$ where we use the usual convention for the multi-indexes $I=(i_1,\dots,i_d),
J=(j_1,\dots,j_d)\in \mathbb N^d$. Then an easy computation yields that 
\begin{enumerate}
\item $k_m(p)= \frac{-1}{m+1}K^1_m(p,p) = -\sum_{|I|=m+1}K_m^{e_1,I-e_1}\frac{p^I}{I!}$ where $e_1=(1,0,\dots,0)$
\item $dK_m = 0$ implies that $K_m^{I,J}= K_m^{L,N}$ if $|I|+|J|=|L|+|N|$
\end{enumerate}
which implies that $dk_m(p_1,p_2)= K_m(p_1,p_2)$.

\end{proof}

\section{Comparison with deformation quantization}\label{kontform}

In this section we make precise the statement that the generating function 
may be seen as the semi-classical approximation of the Kontsevich deformation 
formula. Namely Kontsevich gave in \cite{kontsevich1997} an explicit formula for
the associative deformation of the usual product of function on $\mathbb R^d$ into the direction
of a Poisson structure $\alpha$,
$$f*g = fg +\sum_{n\geq 1} \frac{h^n}{n!}\sum_{\Gamma\in G_{n,2}}W_\Gamma B_\Gamma(f,g)$$
where $W_\Gamma$ are the weights and $B_\Gamma$ the bidifferential operators
introduced in Section \ref{sec:Konttrees}.

\begin{Definition}
Consider a graph $\Gamma$ in $C_{n,2}$, the set of connected graphs  of type  $(n,2)$. 
We denote by $n_{\Gamma}:=|E_\Gamma^a|$ the number of aerial edges and $e_{\Gamma}:=|E_\Gamma^g|$ the number of ground edges.
\end{Definition}

In order to introduce the number of loops in a connected graph
let us make the following remark.
If $\Gamma$ is a connected graph of type  $(n,2)$ then $\Gamma$ must at least have $n-1$ aerial edges.
Which means that $n-1 \leq n_{\Gamma}$. On the other
hand we have $n_{\Gamma}+e_{\Gamma}=2n$ This implies 
that for connected Kontsevich graphs the number $n-e_{\Gamma}
+1$ is always positive or zero.

\begin{Definition}
For a connected graph  of type $(n,2)$ we call the number
$n-e_{\Gamma}+1$ the number of loops of the graph and
we denote it by $b_{\Gamma}$. We denote by $B^l_n$ the set of 
connected graphs of type $(n,2)$  with $l$ loops and we set $B^l= \cup_{n=1}^\infty B^l_n$.
It is easy to see that $B_n^0$ are exactly the Kontsevich trees $T_{n,2}$.
\end{Definition}

The following Lemma shows that the star-product can be considered
as a suitable exponentiation of a deformation of the Poisson structure.

\begin{Lemma}[\bf{Exponential formula}]
Let $f,g\in C^{\infty}(M)$. The star-product could be expressed as
$$f*g(x)=\exp\bigg(\frac{1}{h}D\Big(h \partial_{x'},h\partial_{x''},x\Big)\bigg)
f(x')g(x'')\bigg|_{x'=x''=x}$$ where $D(p_1,p_2,x)=\sum_{j=0}^{\infty}h^{j}D^{j}(p_1,p_2,x)$
and $ D^{j}(p_{1},p_{2},x)=\sum_{\Gamma\in B^{j}}\frac{W_{\Gamma}}{|\Gamma|!}\hat{B}_{\Gamma}(p_{1},p_{2},x)$.

\end{Lemma}

\begin{proof}

By definition of the star-product, the definition of the $\hat{B}$ and using
Lemma \ref{lem:fact} of Section \ref{sec:Konttrees} we can do the following computation,

\begin{eqnarray*}
I   &= &  f*g(x)\\
    &= &  \big\{1+\sum_{n=1}^{\infty} \frac{h^{n}}{n!} \sum_{\Gamma \in G_{n,2}} W_{\Gamma}\hat{B}_{\Gamma}(\partial_{x'},\partial_{x''},x)\big\} f(x')g(x'')\bigg|_{x'=x''=x} \\
    &=&  \big\{1+\sum_{n=1}^{\infty} \frac{h^{n}}{n!}\sum_{\substack{\Gamma \in G_{n,2} \\ D(\Gamma)=(\Gamma_{1},\dots,\Gamma_{k})}} (W_{\Gamma_{1}}\hat{B}_{\Gamma_{1}}) \dots (W_{\Gamma_{k}}\hat{B}_{\Gamma_{k}})\big\}  f(x')g(x'')\bigg|_{x'=x''=x}   \\
    &=&\big\{1+\sum_{n=1}^{\infty} \frac{h^{n}}{n!}\sum_{k=1}^{n}\sum_{\substack{n_{1},\dots,n_{k}\in \mathbb{N}\backslash \{0\} \\ n_{1}+\dots+n_{k}=n}} \sum_{\Gamma \in G_{n,2}(n_{1},\dots,n_{k})} (W_{\Gamma_{1}}\hat{B}_{\Gamma_{1}}) \dots (W_{\Gamma_{k}}\hat{B}_{\Gamma_{k}})\big\} f(x')g(x'')\bigg|_{x'=x''=x} \\
   &=&\big\{1+\sum_{n=1}^{\infty} \frac{h^{n}}{n!}\sum_{k=1}^{n}\sum_{\substack{n_{1},\dots,n_{k}\in \mathbb{N}\backslash\{0\} \\ n_{1}+\dots+n_{k}=n}} (\frac{n!}{k!n_{1}!\dots n_{k}!}) \\
   & & \quad\quad\quad  \sum_{(\Gamma_{1},\dots,\Gamma_{k}) \in C_{n_{1},2}\times \dots \times C_{n_{k},2}} (W_{\Gamma_{1}}\hat{B}_{\Gamma_{1}}) \dots (W_{\Gamma_{k}}\hat{B}_{\Gamma_{k}})\big\} f(x')g(x'')\bigg|_{x'=x''=x} \\
   &=&\big\{1+\sum_{k=1}^{\infty}\frac{1}{k!}\sum_{\substack{(\Gamma_{1},\dots,\Gamma_{k}) \in C_{n_{1},2}\times \dots \times C_{n_{k},2} \\n_{1},\dots,n_{k}\in \mathbb{N}\backslash \{0\} }} (h^{n_{1}}\frac{W_{\Gamma_{1}}}{n_1!}\hat{B}_{\Gamma_{1}}) \dots (h^{n_{k}}\frac{W_{\Gamma_{k}}}{n_k!}\hat{B}_{\Gamma_{k}})\big\} f(x')g(x'')\bigg|_{x'=x''=x} \\
   &=&\big\{1+\sum_{k=1}^{\infty}\frac{1}{k!}(\sum_{n=1}^{\infty}h^{n}\sum_{\Gamma \in C_{n,2}} \frac{W_{\Gamma}}{n!}\hat{B}_{\Gamma})^{k}\big\}f(x')g(x'')\bigg|_{x'=x''=x} \\
   &=&\exp\big\{\frac{1}{h}\sum_{n=1}^{\infty}h^{n+1}\sum_{\Gamma \in C_{n,2}} \frac{W_{\Gamma}}{n!}\hat{B}_{\Gamma}(\partial_{x'},\partial_{x''},x)\big) \big\}f(x')g(x'')\bigg|_{x'=x''=x}
\end{eqnarray*}

Remarking that $\hat{B}_{\Gamma}(\partial_{x'},\partial_{x''},x)
=\frac{1}{h^{e_{\Gamma}}}\hat{B}_{\Gamma}(h\partial_{x'},h\partial_{x''},x)$, we can conclude that

\begin{eqnarray*}
f*g(x) & =& \exp\big\{\frac{1}{h}\sum_{n=1}^{\infty}\sum_{\Gamma \in C_{n,2}} h^{n+1-e_{\Gamma}}\frac{W_{\Gamma}}{|\Gamma|!}\hat{B}_{\Gamma}(h\partial_{x'},h\partial_{x''},x)\big) \big\} \\
       &= &\exp\big\{\frac{1}{h}\sum_{j=0}^{\infty} h^{j}D^{j}(h\partial_{x'},h\partial_{x''},x)\big) \big\}.
\end{eqnarray*}

\end{proof}

The semi-classical part of the deformation formula is $$\frac{1}{h}D^0(h p_{1},h p_{2},x).$$
It is easy to see that 
$$x(p_{1}+p_{2})+\frac{1}{h}D^0(h p_{1},h p_{2},x)$$
is exactly the formal symplectic groupoid generating function. It is in this sense
that one can consider the generating function as a semi-classical approximation of the deformation
formula.

\par
We give now a quick but non rigorous proof of the existence part of Theorem \ref{th:main}. We use
 the technique of saddle point approximation (over non really-well defined integrals). The following computations are then
 by no way a replacement of the rigorous and more technical argument developed in the previous sections.

First notice that as consequence of the exponential formula of the previous Lemma we have that
$$e^{\frac{i}{\hbar}p_1x}*e^{\frac{i}{\hbar}p_2x}= e^{\frac{i}{h}(\sum_{j\geq 0}(\frac{\hbar}{i})^jD^j(p_1,p_2,x))}.$$
We have replaced in the above identity the previously used formal parameter $h$ by $\frac{\hbar}{i}$
for better agreement with the notations in quantum mechanics.
 Moreover we have absorbed the term $x(p_1+p_2)$ into $D^0$. We keep
using this convention through the following computation.
        Let us compute both sides of
        $$\underbrace{(e^{\frac{i}{\hbar}p_{1}x}*e^{\frac{i}{\hbar}p_{2}x})*e^{\frac{i}{\hbar}p_{3}x}}_{A}=
        \underbrace{e^{\frac{i}{\hbar}p_{1}x}*(e^{\frac{i}{\hbar}p_{2}x}e^{\frac{i}{\hbar}p_{3}x})}_{B}$$
        with the help of the asymptotical Fourier transform. We have then,
        $$A=(2\pi \hbar)^ {-d / 2}\int \hat{f}(p_{1},p_{2},p)(e^{\frac{i}{\hbar}px}*e^{\frac{i}{\hbar}p_{3}x})dp$$
        where$\hat{f}(p_{1},p_{2},p)$ is the Fourier transform of
        $$f(p_{1},p_{2},x)=e^{\frac{i}{\hbar}p_{1}x}*e^{\frac{i}{\hbar}p_{2}x}=e^{\frac{i}{\hbar}\sum_{j=0}^{\infty}
        (\frac{\hbar}{i})^{j}D^{j}(p_{1},p_{2},x)}$$
        that is,
        $$\hat{f}(p_{1},p_{2},p)=(2\pi \hbar)^ {-d / 2}\int e^{\frac{i}{\hbar}\big(\sum_{j=0}^{\infty}
        (\frac{\hbar}{i})^{j}D^{j}(p_{1},p_{2},x)-px\big)}dx.$$

        We  use the method of the saddle point approximation to evaluate this integral when ``$\frac{\hbar}{i}$ is very
        small".
        
        First notice that for functions of the form
        $$g_{\frac{\hbar}{i}}(x)=g_{0}(x)+\frac{\hbar}{i} g_{1}(x)+(\frac{\hbar}{i})^{2} g_{2}(x)+\dots$$
        a formal application of the implicit function theorem to
        $F(\frac{\hbar}{i},x)=\nabla g_{\frac{\hbar}{i}}(x)$ tells us that

        \begin{enumerate}
                \item{
                $\exists \bar{x}:I\rightarrow \mathbb{R}^{n}$ where $I$
                is a interval around zero so that $\bar{x}(\frac{\hbar}{i})$ is
                an extremal point of $g_{\frac{\hbar}{i}}(x)$ if $\bar{x}(0)=
                \bar{x}$ is an extremal point of $g_{0}(x)$.
                }
                \item{
                $\bar{x}(\frac{\hbar}{i})=\bar{x}-g_{0}^{-1''}(\bar{x})g_{1}'(\bar{x})\frac{\hbar}{i}
                +\mathcal{O}((\frac{\hbar}{i})^{2})$
                }
        \end{enumerate}

        Second, notice that we have the following asymptotical expansion
        \vspace{10pt} 
        $$g_{\frac{\hbar}{i}}(\bar{x}(\frac{\hbar}{i}))=g_{\frac{\hbar}{i}}((\bar{x}-
        g_{0}^{-1''}(\bar{x})g_{1}'(\bar{x})\frac{\hbar}{i}+\mathcal{O}((\frac{\hbar}{i})^{2})),$$
        around $\bar{x}$ we get
        $$g_{\frac{\hbar}{i}}(\bar{x}(\frac{\hbar}{i}))=g_{0}(\bar{x})+\frac{\hbar}{i} g_{1}(\bar{x})
        +\mathcal{O}((\frac{\hbar}{i})^{2}).$$

        Now if we apply the method of the stationary phase to
        $$I=\int e^{\frac{i}{\hbar}g_{\frac{\hbar}{i}}(x)}dx$$
        we find
        $$I\approx c(\bar{x},\frac{\hbar}{i})e^{\frac{i}{\hbar}(g_{0}(\bar{x})+
        \frac{\hbar}{i} g_{1}(\bar{x}))}$$
        where $\bar{x}$ is the extremal point of $g_{0}$.

        Let us come back to the computation  of $A$.
        With the preceding remarks in mind the computation of $\hat{f}(p_{1},p_{2},p)$
        leads, through the application of the stationary phase method, to
        $$
        \hat{f}(p_{1},p_{2},p)\approx c(p_{1},p_{2},\bar{x},\frac{\hbar}{i})
        e^{\frac{i}{\hbar}(D^0(p_{1},p_{2},\bar{x})-p\bar{x}+\frac{\hbar}{i} 
        D^1(p_{1},p_{2},\bar{x}))}
        $$

        where $c$ is a certain function of $p_{1},p_{2}$, $\bar{x}$ 
        and $\frac{\hbar}{i}$ and where $\bar{x}$ is a critical point of
        $$D^0(p_{1},p_{2},x)-px \quad (i.e.,\quad\nabla_x D^0(p_{1},p_{2},\bar{x})=p).$$ Then
        $$A \approx \int c(p_{1},p_{2},\bar{x},\frac{\hbar}{i})e^{\frac{i}{\hbar}
        (D^0(p_{1},p_{2},\bar{x})-p\bar{x}+\frac{\hbar}{i} D^1(p_{1},p_{2},\bar{x})+
        \sum_{j=0}^{\infty}(\frac{\hbar}{i})^{j}D^{j}(p,p_{3},x))}dp.$$

        Using the same method as above again we obtain
        $$A \approx \tilde{C}(p_{1},p_{2},\bar{p},\bar{x},x,\frac{\hbar}{i})
        e^{\frac{i}{\hbar}(D^0(p_{1},p_{2},\bar{x})-\bar{p}\bar{x}
        +D^0(\bar{p},p_{3},x))} e^{-D^1(\bar{p},p_{3},x)
        -D^1(p_{1},p_{2},\bar{x})}$$
        where $\bar{x}$ is determined by $\nabla_x D^0(p_{1},p_{2},\bar{x})=\bar{p}$
        as above and $\bar{p}$ by $\nabla_{P_1}D^0(\bar{p},p_{3},x)=\bar{x}$.

        Namely, $$\frac{d}{dp}\big[D^0(p_{1},p_{2},\bar{x})-\bar{p}\bar{x}+D^0(\bar{p},p_{3},x)\big] = 0 $$ gives 
         $$ \nabla_x D^0(p_{1},p_{2},\bar{x})\frac{d \bar{x}} {dp}-
         \nabla_{x}D^0(p_{1},p_{2},\bar{x})\frac{d \bar{x}} {dp}-
         \bar{x}+\nabla_{p_1}D^0(\bar{p},p_{3},x) =  0 .$$

        By the same kind of computation we approximate $B$ for $\frac{\hbar}{i}$ ``small enough",
        $$B \approx \tilde{C}(p_{2},p_{3},\tilde{p},\tilde{x},x,\frac{\hbar}{i})
        e^{\frac{i}{\hbar}(D^0(p_{2},p_{3},\tilde{x})-\tilde{p}\tilde{x}
        +D^0(p_{1},\tilde{p},x))}
        e^{-D^1(p_{2},p_{3},\tilde{x})-D^1(p_{1},\tilde{p},x)}
        $$
        with $\tilde{x}$ and $\tilde{p}$ determined by $\nabla_x D^0(p_{2},p_{3},\tilde{x})=\tilde{p}$
        and $\nabla_{p_2}D^0(p_{1},\tilde{p},x)=\tilde{x}$.

        \vspace{10pt}
        Equating $A$ and $B$ we then get that $D^0(p_1,p_2,x)$ satisfies the SGA equation.

\end{document}